\documentclass{article}

\usepackage{arxiv}

\usepackage[utf8]{inputenc}
\usepackage[T1]{fontenc}

\usepackage{amsmath,amssymb,amsbsy}
\usepackage{amsthm}
\usepackage{newtxmath}
\usepackage{bm}
\usepackage{hyperref}
\usepackage{url}
\usepackage{booktabs}
\usepackage{amsfonts}
\usepackage{nicefrac}
\usepackage{microtype}
\usepackage{graphicx}
\usepackage{natbib}
\usepackage{doi}

\usepackage{xcolor}
\usepackage{enumitem}

 \bibpunct[, ]{(}{)}{,}{a}{}{,}%
\theoremstyle{plain}
\newtheorem{theorem}{Theorem}
\newtheorem{proposition}{Proposition}
\newtheorem{assumption}{Assumption}
\theoremstyle{definition}
\newtheorem{remark}{Remark}

\newcommand{\calT}{\mathcal{T}}
\newcommand{\calX}{\mathcal{X}}
\newcommand{\calD}{\mathcal{D}}

\newcommand{\E}{\mathbb{E}}
\newcommand{\R}{\mathbb{R}}

\newcommand{\bmu}{\boldsymbol{\mu}}
\newcommand{\bpsi}{\boldsymbol{\psi}}
\DeclareMathOperator{\conv}{conv}
\DeclareMathOperator{\dom}{dom}
\DeclareMathOperator{\logit}{logit}

\title{\upshape\bfseries Decision-Focused Learning for Multistage Stochastic Programs,\\ with Application to Power Network Repair Scheduling}

\author{\href{https://orcid.org/0009-0004-0547-1941}{\includegraphics[scale=0.06]{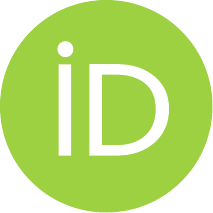}\hspace{1mm}David Krame Kadurha}, \href{https://orcid.org/0000-0001-6585-4785}{\includegraphics[scale=0.06]{orcid.pdf}\hspace{1mm}Alessandro Zocca}, \href{https://orcid.org/0000-0003-1124-8821}{\includegraphics[scale=0.06]{orcid.pdf}\hspace{1mm}Sandjai Bhulai} \\
	Department of Mathematics, Vrije Universiteit Amsterdam, Amsterdam, The Netherlands \\
	\texttt{\{d.kramekadurha, a.zocca, s.bhulai\}@vu.nl}
}

\renewcommand{\shorttitle}{DFL for MSPs, with Application to Power Network Repair Scheduling}
\renewcommand{\headeright}{}
\renewcommand{\undertitle}{}

\definecolor{linkblue}{RGB}{0,60,250}
\hypersetup{
colorlinks=true,
allcolors=linkblue,
pdftitle={Decision-Focused Learning for Multistage Stochastic Programs, with Application to Power Network Repair Scheduling},
pdfauthor={David Krame Kadurha, Alessandro Zocca, Sandjai Bhulai},
pdfkeywords={multistage stochastic programming, decision-focused learning, Fenchel--Young loss, perturbed optimizers, power network repair scheduling},
}

\begin{document}
\maketitle

\begin{abstract}
Multistage stochastic programs (MSPs) with discrete decisions are hard to solve due to the curse of dimensionality, non-convex value functions, and scenario trees that grow exponentially with the horizon. We develop a decision-focused learning (DFL) framework for MSPs in which a neural network predicts scores for the part of the decision that drives the state transition. Under structural conditions that make this sub-decision enter the stage objective linearly, the network is trained on labels from offline single-scenario problems, using the Fenchel--Young loss with perturbed optimizers, which need only an optimization oracle. We apply the framework to repair scheduling in power networks under stochastic, time-varying component failures, with and without a per-stage repair budget. On a synthetic 5-bus network and on the IEEE 14-bus network at two repair-cost scales, the learned and baseline policies are compared on common evaluation paths. Without the budget, DFL is the cheapest policy on the 5-bus network; on the 14-bus network it does not outperform a near-sighted model-predictive control (MPC) baseline, but it beats repairing every failed component at the higher scale. Under a binding budget, DFL is the cheapest policy evaluated on both networks and at both scales, 14\% to 23\% below the next-best policy. On the 5-bus network, the mean cost of DFL is 1.7\% (without the budget) and 1.8\% (with it) above the exact optimum computed by dynamic programming. Within several days of computation, sample average approximation yields costlier policies than DFL and the repair-everything heuristic, and stochastic dual dynamic integer programming remains far from convergence. DFL computes its decisions about as fast as the trivial heuristics and orders of magnitude faster than the MPC baseline. Since these conditions are not specific to power networks, the framework extends to other MSPs that share this structure. Code is available at \url{https://github.com/DavidKrame/DFL_for_MSP__Application_to_power_network_repair}.
\end{abstract}


\section{Introduction}
\label{sec:introduction}

Multistage stochastic programs (MSPs) model decisions that adapt to information revealed over time and arise in, among other fields, inventory management, energy dispatch, financial optimization, supply chain planning, and, as studied here, power grid operations under uncertainty. Their exact solution is hindered by three related difficulties: value functions must be represented over a state space whose size grows exponentially with its dimension; integer decisions can make these functions non-convex; and a scenario tree with $q$ branches per stage has $q^T$ leaves, a number exponential in the horizon $T$. SDDP \citep{pereira1989optimal,pereira1991multi} relies on convex stage problems and does not directly accommodate integer recourse. Its integer extension SDDiP \citep{zou2019sdip} converges finitely for binary state variables, so continuous states must be approximated by binarization, and its worst-case effort grows exponentially with the number of binary state variables. SAA \citep{kleywegt2002saa,shapiro2009lectures} produces deterministic mixed-integer programs that grow with the scenario tree.

We ask whether the part of the continuation value that the current decision needs can be learned from offline data, without reconstructing the value function or building a scenario tree. Building on decision-focused learning (DFL) \citep{elmachtoub2022smart,mandi2024decision}, under the structural conditions stated in Section~\ref{sec:dfl_general}, we parameterize the continuation value linearly in a learnable sub-decision and train the resulting optimization layer against decision quality. Training labels come from offline, single-scenario perfect-information MIPs, as is common in optimization-augmented machine learning \citep{Maximilian_DFL_Survey}. These MIPs grow only linearly with the horizon, and deployment avoids both backward recursion and scenario trees.

We instantiate the framework for repair scheduling in power networks under time-varying stochastic component failures, motivated by severe-weather risk. At each stage, the operator dispatches power under DC Optimal Power Flow constraints and schedules repairs. Structural analysis (Section~\ref{subsec:decoupling}) yields a closed-form threshold repair rule at each stage, driven by learned continuation-value coefficients that depend on the state and the stage; a per-stage budget instead couples the repairs of a stage through a knapsack constraint.

Our contributions are threefold. First, we present a general DFL framework for MSPs, identify conditions supporting a linear continuation-value parameterization, and provide a training and deployment procedure. Second, we formulate power-network repair scheduling under stochastic component failures as a worked example and identify the decomposition and decoupling properties that make the framework applicable. With the budget, the knapsack surrogate is not separable, so training uses the general perturbed-optimizer estimator of Section~\ref{subsec:fy}. Third, experiments on synthetic 5-bus and IEEE 14-bus networks show that DFL is the cheapest policy evaluated under a binding per-stage budget, that its advantage over repairing every failed component grows with the repair-cost scale, and that it deploys at the cost of a simple heuristic.

The paper is organized as follows. Section~\ref{sec:literature} reviews related work, followed by Section~\ref{sec:dfl_general}, which develops the general DFL framework for MSPs. Then, in Section~\ref{sec:formulation}, we formulate the repair scheduling problem and identify its structural properties. Section~\ref{sec:experiments} presents computational experiments, and finally, Section~\ref{sec:conclusion} concludes. The appendix collects the proofs and the Fenchel--Young background (Section~\ref{ec:fy}), the learning implementation (Section~\ref{ec:impl}), the full model and failure process (Section~\ref{ec:model}), the specifications of the comparison methods (Section~\ref{ec:comparators}), and how each method used the computing resources (Section~\ref{ec:compute}).

\section{Literature Review}
\label{sec:literature}

\paragraph{Multistage stochastic programming.} The MSP framework, including the value functions on which the solution methods below rest, is developed in \citet{shapiro2009lectures}. The dominant exact approaches represent the future explicitly, either through cuts built by backward recursion, as in SDDP \citep{pereira1989optimal,pereira1991multi} and SDDiP \citep{zou2019sdip}, or through an explicit scenario tree, as in SAA \citep{kleywegt2002saa}. Their cost is driven by the number of cuts or scenario-tree nodes required before that representation is accurate enough, which with discrete decisions runs into the difficulties summarized in Section~\ref{sec:introduction}.

\paragraph{Decision-focused learning.} Our learning machinery comes from decision-focused learning, whose premise is the predict-then-optimize paradigm of \citet{elmachtoub2022smart}: rather than minimizing prediction error, one trains directly against the quality of the resulting decision. The survey of \citet{mandi2024decision} gives a comprehensive account. The central technical obstacle is differentiating through a combinatorial optimization layer, and it has been addressed through Quadratic Program (QP) differentiation via the implicit function theorem \citep{amos2017optnet}, the Fenchel--Young loss framework \citep{blondel2020learning}, and perturbed optimizers \citep{berthet2020learning}, the latter two connected through Fenchel duality in \citet{berthet2020learning} and set in a general probabilistic framework by \citet{dalle2022learning}. These tools, however, are almost always deployed in a \emph{single-stage} setting, where one prediction feeds one downstream optimization, and only rarely in a \emph{two-stage} one; in multistage settings they have, to our knowledge, appeared only in problem-specific forms \citep{Axel_VRPTW,Maximilian_DFL_EMS}. The closest step toward a general multistage method is decision-focused forecasting \citep{EgonMiguel2025}, which forecasts the exogenous parameters entering a multistage program and backpropagates decision regret through the entire solve; it is restricted, however, to convex stage problems with a differentiable solution map and does not treat integer or combinatorial stage problems.

\paragraph{Power network optimization.} Our worked example builds on the DC optimal power flow model \citep[reviewed in][]{wood2013power}, written in the network-agnostic operator form of \citet{zhou2021opf}, which we adopt and extend to a multistage setting with repair decisions. Service restoration after outages has been studied mostly for distribution grids \citep{eisenkraemer2026systematic,goda2025restoration}. The influence of adverse weather on component failure rates is examined in \citet{panteli2015influence}, and time-varying failure probability models under severe weather are developed in \citet{alvehag2011reliability} and \citet{wang2002reliability}. To schedule repairs against these failures, that literature relies largely on heuristics or on stochastic and robust MILP formulations, whose cost grows quickly with the number of components and the planning horizon. We instead cast repair scheduling under such stochastic, time-varying failures as an instance of our multistage DFL framework, and we benchmark it against an exact multistage method by running SDDiP through the \texttt{SDDP.jl} implementation of \citet{dowson2021sddpjl}.

\paragraph{Positioning.} The closest works are \citet{Axel_VRPTW} on dynamic vehicle routing and \citet{Maximilian_DFL_EMS} on ambulance dispatching and redeployment. In both, a learned model predicts the parameters of an optimization layer designed for the application (prizes in a prize-collecting routing problem, edge weights in a bipartite matching) and is trained with a perturbed loss to imitate perfect-information solutions. We use the same training scheme, but our layer comes from the MSP itself. For any MSP that meets two structural conditions (Section~\ref{sec:dfl_general}), it is the stage problem in which the continuation value is replaced by a linear function of the sub-decision driving the state transition, with learned coefficients. When the optimal cost of the remaining decisions does not depend on that sub-decision, as in our application, the layer is a linear-objective oracle, so the Fenchel--Young loss and perturbed-optimizer gradients apply without making the solver differentiable, even when the stage problems are mixed-integer. Unlike the classical MSP methods, the approach uses no backward pass or scenario tree, and at deployment it solves one stage problem per stage. Our contribution is this general construction, the conditions under which it applies, and a worked application; the Fenchel--Young and perturbation results we rely on are established \citep{blondel2020learning,berthet2020learning,dalle2022learning}, and we collect them with their proofs in Section~\ref{ec:fy}.

\section{A DFL Framework for Multistage Stochastic Programs}
\label{sec:dfl_general}

\subsection{Multistage stochastic programming}
\label{subsec:msp}

Let $\calT = \{1, \ldots, T\}$ be a finite planning horizon. At each stage $t \in \calT$ the system occupies a state $x_t \in \calX$, an exogenous random variable $\xi_t \in \Xi$ is revealed, the controller selects a decision $\delta_t \in \calD_t(x_t, \xi_t)$, pays stage cost $c_t(x_t, \delta_t, \xi_t)$, and the state evolves as $x_{t+1} = f_t(x_t, \delta_t, \xi_t)$. A \emph{policy} $\pi = \{\pi_t\}_{t=1}^T$ is a collection of decision rules $\pi_t$, each selecting a feasible decision $\delta_t = \pi_t(x_t,\xi_{1:t}) \in \calD_t(x_t,\xi_t)$ from the information available at stage $t$. The
\emph{non-anticipativity} requirement is that this choice depend only on the realizations $\xi_{1:t}$ observed up to stage $t$, never on future realizations $\xi_{t+1:\!T}$; equivalently, each decision $\delta_t$ is $\sigma(\xi_{1:t})$-measurable. The MSP minimizes total expected cost over all non-anticipative policies:
\begin{equation}
  \min_\pi \; \E\!\left[\sum_{t=1}^T c_t(x_t,\pi_t(x_t,\xi_{1:t}),\xi_t)\right].
  \label{eq:msp}
\end{equation}

We assume throughout that~\eqref{eq:msp} satisfies \emph{complete recourse}: for every $t \in \calT$, $x_t \in \calX$, $\xi_t \in \Xi$, the feasible set $\calD_t(x_t,\xi_t)$ is non-empty. We also assume that the shocks $\xi_1,\ldots,\xi_T$ are stagewise independent, so that $x_t$ is a Markov state, and that the minima below are attained; both conditions hold in our application. Under these conditions and standard measurability and integrability requirements \citep{shapiro2009lectures}, the value functions $\mathcal{V}_t$ are well defined and the optimal value of~\eqref{eq:msp} satisfies the Bellman equations with $\mathcal{V}_{T+1} \equiv 0$:
\begin{equation}
  \mathcal{V}_t(x_t) = \E_{\xi_t}\!\left[
    \min_{\delta_t \in \calD_t(x_t,\xi_t)}
    \bigl\{ c_t(x_t,\delta_t,\xi_t) + \mathcal{V}_{t+1}(f_t(x_t,\delta_t,\xi_t)) \bigr\}
  \right].
  \label{eq:bellman}
\end{equation}
An optimal policy is then any non-anticipative rule $\pi_t^*$ that selects, at every stage and state, a minimizer of the inner problem in~\eqref{eq:bellman}. We will refer to $\mathcal{V}_{t+1}(x_{t+1}) = \mathcal{V}_{t+1}(f_t(x_t,\delta_t,\xi_t))$ as the \emph{continuation value} at state $x_{t+1}$; it equals the expected optimal cost from stage $t+1$ onwards, starting from state $x_{t+1}$.

\subsection{Linear continuation-value parameterization and the DFL surrogate}
\label{subsec:surrogate}

Instead of reconstructing $\mathcal{V}$, our approach identifies, within the problem structure, a low-dimensional parameterization of the continuation value that is \emph{linear in the decisions of interest} and that can be learned from offline data.

Assume the decision $\delta_t$ decomposes as $\delta_t = (\delta_t^0, u_t)$ into a \emph{learnable sub-decision} $u_t$ and a complementary part $\delta_t^0$.   The DFL framework requires two structural conditions on this decomposition.

\begin{assumption}[State-transition separability]
\label{assump:state_dep} The state transition depends on the full decision $(\delta_t^0, u_t)$ only through $u_t$:
\[ f_t(x_t, (\delta_t^0, u_t), \xi_t) = f_t(x_t, u_t, \xi_t).
\]
\end{assumption}

\begin{assumption}[Linear stage cost in $u_t$]
\label{assump:linear_cost} The stage cost decomposes as
\[ c_t(x_t, \delta_t, \xi_t) = c_t^0(x_t, \delta_t^0, \xi_t) + \beta_t^\top u_t
\] for a fixed cost vector $\beta_t \in \R^{|u_t|}$, where $c_t^0$ does not depend on $u_t$.
\end{assumption}

Denote the marginal feasible set for $u_t$ by
\begin{equation}
  \mathcal{U}_t(x_t,\xi_t)
  \;:=\;
  \bigl\{u_t : \exists\,\delta_t^0 \text{ s.t. } (\delta_t^0,u_t) \in \calD_t(x_t,\xi_t)\bigr\}.
  \label{eq:U_def}
\end{equation}
Under Assumptions~\ref{assump:state_dep} and~\ref{assump:linear_cost}, the Bellman equation~\eqref{eq:bellman} factors as
\begin{equation}
  \mathcal{V}_t(x_t)
  = \E_{\xi_t}\!\left[
    \min_{u_t \in \mathcal{U}_t(x_t,\xi_t)}
      \Bigl\{
        Q_t^0(x_t,u_t,\xi_t)
        + \beta_t^\top u_t
        + \mathcal{V}_{t+1}(f_t(x_t,u_t,\xi_t))
      \Bigr\}
  \right],
  \label{eq:bellman_decoupled}
\end{equation}
where
\begin{equation}
  Q_t^0(x_t,u_t,\xi_t)
  \;:=\;
  \min_{\substack{\delta_t^0:\;(\delta_t^0,u_t)\in\mathcal{D}_t(x_t,\xi_t)}}
  c_t^0(x_t,\delta_t^0,\xi_t)
  \label{eq:inner_cost}
\end{equation}
is the optimal inner cost for a fixed $u_t$. The feasible set in~\eqref{eq:inner_cost} may depend on $u_t$, in which case $Q_t^0$ is itself a function of $u_t$ and remains inside the outer minimization, while by Assumption~\ref{assump:state_dep} the continuation value depends on the decision only through $u_t$.

The factorization~\eqref{eq:bellman_decoupled}--\eqref{eq:inner_cost} is exact under the two assumptions. The framework then approximates the continuation value by a function that is affine in $u_t$:
\begin{equation}
  \mathcal{V}_{t+1}(f_t(x_t, u_t, \xi_t)) \approx
  h_t(x_t, \xi_t) - \bmu_t^\top u_t~,
  \label{eq:linear_approx}
\end{equation}
where $\bmu_t \in \R^{|u_t|}$ are the \emph{continuation-value coefficients}, and $h_t(x_t, \xi_t)$ is an additive baseline independent of $u_t$ (it does not affect the surrogate optimization). Substituting into \eqref{eq:bellman_decoupled} and dropping the constant $h_t(x_t,\xi_t)$, the stage surrogate for $u_t$, which is the optimization layer of our framework, is
\begin{equation}
  \min_{u_t \in \mathcal{U}_t(x_t,\xi_t)}
  \bigl\{Q_t^0(x_t,u_t,\xi_t) + (\beta_t - \bmu_t)^\top u_t\bigr\}.
  \label{eq:surrogate_u}
\end{equation}
This is a (possibly combinatorial) optimization problem parameterized by $\bmu_t$. A neural network $\pi_{\bpsi}$ maps the observable state to $\bmu_t$; the network is trained so that solving~\eqref{eq:surrogate_u} with the predicted $\hat{\bmu}_t = \pi_{\bpsi}(x_t, \xi_t)$ reproduces the perfect-information labels $u_t^*$, obtained by solving the horizon-wide deterministic problem on a single sampled scenario (Section~\ref{subsec:train_deploy}).

\begin{remark}[Linear-objective oracle]
\label{rem:fy_connection} The Fenchel--Young machinery of Section~\ref{subsec:fy} needs only an oracle that maximizes a \emph{linear} objective over a compact feasible set, which may be binary, discrete, continuous, or mixed. The surrogate~\eqref{eq:surrogate_u} is such an oracle when $Q_t^0$ is independent of $u_t$, or affine in it with its coefficient absorbed into the score: in the independent case it reads $\arg\max_{u_t \in \mathcal{U}_t} \eta_t^\top u_t$ with \emph{effective score} $\eta_t := \bmu_t - \beta_t$.
\end{remark}

\begin{remark}
\label{rem:assump1} Together with the condition on $Q_t^0$ in Remark~\ref{rem:fy_connection}, Assumptions~\ref{assump:state_dep} and~\ref{assump:linear_cost} are sufficient, not necessary, for a linear-objective surrogate. Assumption~\ref{assump:state_dep} makes the continuation value a function of $u_t$ alone, and Assumption~\ref{assump:linear_cost} makes the known cost linear in $u_t$, so the affine correction~\eqref{eq:linear_approx} enters the surrogate linearly. DFL has also been applied where a linear-objective surrogate of the downstream problem is identified by other means, as in the structured approximation scheme of \citet{parmentier2021}, its multistage instantiation for dynamic vehicle routing in \citet{Axel_VRPTW}, and the variants collected in the survey of \citet{mandi2024decision}. The assumptions say nothing about how \emph{accurate} the affine form~\eqref{eq:linear_approx} is, a question we return to in Section~\ref{sec:conclusion}; the present paper stays within this linear case and leaves non-linear surrogates to future work.
\end{remark}

\subsection{Learning with the Fenchel--Young loss and perturbed optimizers}
\label{subsec:fy}

Let $\mathcal{Y} \subseteq \R^m$ be the compact feasible set of the oracle, $\mathcal{C} = \conv(\mathcal{Y})$ its convex hull, and $\Omega: \mathcal{C} \to \R$ a convex regularizer with conjugate $\Omega^*(\eta) = \max_{u \in \mathcal{C}}\{\eta^\top u - \Omega(u)\}$. The \emph{Fenchel--Young (FY) loss} of \citet{blondel2020learning} is
\begin{equation}
  L_\Omega(\eta;\, u) := \Omega^*(\eta) + \Omega(u) - \eta^\top u
  \;\geq\; 0,
  \label{eq:fy}
\end{equation}
and, when $\Omega$ is strongly convex, its gradient is $\nabla_{\eta} L_\Omega(\eta; u) = \hat{u}_\Omega(\eta) - u$, where $\hat{u}_\Omega(\eta)$ is the maximizer in $\Omega^*(\eta)$ (Proposition~\ref{prop:fy_props} in Section~\ref{ec:fy}). Computing $\hat{u}_\Omega$ over the convex hull is itself a hard problem. \citet{berthet2020learning} instead use the \emph{perturbed optimizer}: fix $\varepsilon > 0$, let $Z$ be a random vector whose distribution $\nu$ has a positive, differentiable density on $\R^m$, a finite first moment, and $\E[Z] = \mathbf{0}$, and define
\begin{equation}
  \hat{u}_\varepsilon(\eta) := \E_{Z\sim\nu}\!\left[
    \arg\max_{u \in \mathcal{Y}} (\eta + \varepsilon Z)^\top u\right].
  \label{eq:perturbed}
\end{equation}
One estimates it without bias by Monte Carlo: sampling $Z_1,\ldots,Z_M \overset{\text{i.i.d.}}{\sim}\nu$, solving $\arg\max_{u \in \mathcal{Y}}(\eta+\varepsilon Z_m)^\top u$ for each $m$, and averaging. Each sample requires a single call to any solver that maximizes a linear objective over $\mathcal{Y}$, with no solver modifications. For the regularizer $\Omega$ induced by $\nu$, $\hat{u}_\varepsilon$ is the regularized decision of $\varepsilon\Omega$, so $\hat{u}_\varepsilon(\eta) - u$ is the exact gradient of the FY loss $L_{\varepsilon\Omega}(\eta; u)$ (Theorem~\ref{thm:equiv} in Section~\ref{ec:fy}; \citealt{berthet2020learning,dalle2022learning}). The FY loss can therefore be minimized without forming $\Omega$, by solving problems with randomly perturbed linear objectives over $\mathcal{Y}$.

\subsection{Training and deployment}
\label{subsec:train_deploy}
Training scenarios may be sampled from a known distribution of $\xi$ or taken directly from historical data.
\paragraph{Training.} For each training scenario $s = 1,\ldots,N_\text{train}$, solve the single-scenario MIP (the deterministic $T$-stage problem with fixed realization $\xi_s = (\xi_s^1,\ldots,\xi_s^T)$) to obtain the \emph{perfect-information labels} $\{u_s^{*,t}\}_{t=1}^T$. Each such MIP sees the whole sampled trajectory $\xi_s$ before deciding, so these labels are optimal in hindsight for that trajectory and are not draws from an optimal non-anticipative policy; their conditional mean need not coincide with the optimal online action. Training non-anticipative policies on such labels follows established practice \citep{Maximilian_DFL_Survey,Axel_VRPTW,Maximilian_DFL_EMS}. Single-scenario MIPs have no scenario branching: their size grows linearly with the horizon $T$ rather than exponentially with the scenario tree, and the offline training cost is discussed in Section~\ref{subsec:time}. The learning goal is to minimize the aggregate FY loss
\begin{equation}
  \mathcal{L}(\bpsi) = \frac{1}{N_\text{train} T}
  \sum_{s,t} L_{\varepsilon\Omega}\!\bigl(
    \eta_s^t;\; u_s^{*,t}\bigr),
  \qquad
  \eta_s^t := \pi_{\bpsi}(x_s^t,\xi_s^t) - \beta_t,
  \label{eq:train_loss}
\end{equation}
by stochastic gradient descent. Here $\eta_s^t$ is the effective score of Remark~\ref{rem:fy_connection}, the same score the oracle uses at deployment, and the gradient with respect to the network output is $\hat{u}_\varepsilon(\eta_s^t) - u_s^{*,t}$ (Section~\ref{subsec:fy}). Training therefore identifies the scores only through the decisions they induce: it fits them so that the oracle reproduces the labels, and the fitted scores do not need to reflect the true continuation values.

The single-scenario MIP may have several optimal solutions, such as symmetric repair schedules, of which the solver returns one; \eqref{eq:train_loss} then targets that particular vector, although the alternatives encode the same scenario cost. Separately, incumbents accepted within the solver's optimality tolerance may be slightly suboptimal, so a label may differ from every exact optimum.

\paragraph{Deployment.} At each stage $t$: (1) observe $x_t$ and $\xi_t$; (2) compute $\bmu_t = \pi_{\bpsi}(x_t, \xi_t)$ via a single forward pass; (3) solve the surrogate problem~\eqref{eq:surrogate_u} for $u_t$; (4) solve the inner problem for $\delta_t^0$. Steps~(3)--(4) together solve one stage of the recourse problem, so deployment proceeds one stage at a time, without a scenario tree or dynamic programming.

\subsection{Optimality of the DFL policy}
\label{subsec:guarantees}

Two sources of error separate the learned policy from an optimal one: the affine form~\eqref{eq:linear_approx}, and the learning step, through the imitation of perfect-information labels, the finite training sample, the capacity of the network class, and the numerical optimization of the loss. The following proposition states when the resulting policy is optimal.

\begin{proposition}[Conditional optimality of the DFL policy]
\label{prop:conditional}
Suppose that, at every stage $t$ and every reachable pair $(x_t,\xi_t)$, there exist $h_t(x_t,\xi_t)$ and $\bmu_t^*(x_t,\xi_t) \in \R^{|u_t|}$ such that $\mathcal{V}_{t+1}(f_t(x_t,u_t,\xi_t)) = h_t(x_t,\xi_t) - \bmu_t^*(x_t,\xi_t)^\top u_t$ for all $u_t \in \mathcal{U}_t(x_t,\xi_t)$, and that $\pi_{\bpsi}(x_t,\xi_t) = \bmu_t^*(x_t,\xi_t)$ at all these pairs. Then, the DFL policy is optimal for~\eqref{eq:msp}.
\end{proposition}

The proof is given in Section~\ref{ec:proofs}. Proposition~\ref{prop:conditional} describes the ideal case that the framework approximates, not a guarantee for the trained network, which is fitted to perfect-information labels (Section~\ref{subsec:train_deploy}). On the 5-bus network, where the optimal expected cost can be computed exactly, the mean cost of the learned policy is 1.7\% above it without a repair budget and 1.8\% above it with one (Section~\ref{sec:experiments}).

\section{Power Network Repair Scheduling}
\label{sec:formulation}

We now instantiate the general framework of Section~\ref{sec:dfl_general} for the power network repair scheduling problem: we adapt a standard DC-OPF to account for component outages, formulate the multistage stochastic problem, and identify the structural properties that make the DFL framework applicable.

\subsection{DC optimal power flow}
\label{subsec:dc_opf}

\paragraph{Starting point.} We adapt the network-agnostic DC-OPF in Equation~(3) of \citet{zhou2021opf}. The directed network has $N$ buses and $N_E$ lines, with bus--line incidence matrix $C\in\{-1,0,1\}^{N\times N_E}$ and diagonal susceptance matrix $B=\operatorname{diag}(b_1,\ldots,b_{N_E})$. It contains $N_G$ generators, with generation $p_g\in\R^{N_G}$, linear cost $c$, and limits $(\underline p_g,\overline p_g)$, and $N_L$ loads with demand $p_l\in\R^{N_L}$; $(\underline f,\overline f)$ are the line limits. Bus-mapping matrices $S_g\in\{0,1\}^{N\times N_G}$ and $S_l\in\{0,1\}^{N\times N_L}$ place each generator and each load at its bus, and under the DC approximation \citep{wood2013power}, voltage angles $\theta\in\R^N$ induce line flows $BC^\top\theta$.

\paragraph{Extension to time-varying availability and load shedding.} Our multistage formulation extends that starting point in two ways. First, component failures change the network: failed generators are taken offline, and failed lines are disconnected, that is, removed from the network so that they carry no power \emph{and} impose no coupling between their endpoint voltage angles. We model this through a time-varying availability vector $a^t = (a_g^t, a_e^t) \in \{0,1\}^{N_G+N_E}$, where $a^t[k] = 1$ if component $k$ is operational at stage $t$ and $a^t[k] = 0$ if it is failed (defined precisely in Section~\ref{subsec:msmip}). Because a line outage is a topology change, we introduce explicit line-flow variables $f^t \in \R^{N_E}$ and impose the DC power-flow law conditionally on line availability. Second, demand may be unserved: we introduce the demand served $d^t \leq p_l^t$ and penalize load shedding at rate $\alpha > 0$. The resulting stage-$t$ DC-OPF is
\begin{equation}
\begin{aligned}
  \min_{p_g^t,\,\theta^t,\,f^t,\,d^t} \quad
  & c^\top p_g^t + \alpha\,\mathbf{1}^\top(p_l^t - d^t) \\[2pt]
  \text{s.t.} \quad
  & \theta^t[\mathrm{ref}] = 0,
    \qquad -\theta^{\max}\mathbf{1} \leq \theta^t \leq \theta^{\max}\mathbf{1}, \\
& C f^t = S_g p_g^t - S_l d^t, \\
  & -(\mathbf{1} - a_e^t) \odot M \;\leq\; f^t - BC^\top\theta^t
    \;\leq\; (\mathbf{1} - a_e^t) \odot M, \\
  & a_e^t \odot \underline{f} \leq f^t \leq a_e^t \odot \overline{f}, \\
  & a_g^t \odot \underline{p}_g \leq p_g^t \leq a_g^t \odot \overline{p}_g, \\
  & \mathbf{0} \leq d^t \leq p_l^t,
\end{aligned}
  \label{eq:dc_opf}
\end{equation}
where $\odot$ denotes element-wise multiplication, $\theta^{\max}$ bounds the bus voltage angles, and $M \in \R^{N_E}_{>0}$ with $M_j := 2\theta^{\max} |b_j|$ is a valid big-$M$ constant, the largest value of $\bigl|b_j(\theta^t[\mathrm{from}(j)] - \theta^t[\mathrm{to}(j)])\bigr|$ over the angle bounds. When $a_e^t[j] = 1$, the two line constraints reduce to the DC power-flow law $f^t[j] = b_j(\theta^t[\mathrm{from}(j)] - \theta^t[\mathrm{to}(j)])$ and the thermal limit $\underline{f}[j] \leq f^t[j] \leq \overline{f}[j]$. When $a_e^t[j] = 0$, the capacity constraint forces $f^t[j] = 0$, while the conditional law relaxes to a bound implied by the angle bounds, so it imposes no relation between the two endpoint angles: the failed line is fully removed from the network and, carrying no flow, drops out of the flow-form balance. Likewise, $a_g^t[i] = 0$ forces $p_g^t[i] = 0$. For fixed $a^t$, \eqref{eq:dc_opf} is therefore the DC-OPF on the energized subnetwork $E^t := \{j : a_e^t[j] = 1\}$; if $E^t$ splits into islands, $\theta^t[\mathrm{ref}] = 0$ pins the angle datum of the island containing the reference bus, and the angles of every other island are determined only up to a common offset, which the angle bounds restrict and which has no effect on the flows, the dispatch, or the cost. When $\underline{p}_g = \mathbf{0}$ (as in all our experiments), the all-off dispatch $p_g^t = d^t = f^t = \mathbf{0}$ (with $\theta^t = \mathbf{0}$) is feasible; hence \eqref{eq:dc_opf} admits complete recourse for any $a^t$.

\subsection{The multistage repair-scheduling problem}
\label{subsec:msmip}

\paragraph{Within-stage chronology.} Each stage $t$ unfolds in a fixed order: the failure vector $\xi^t$ is revealed; the availability updates to $a^t$ (a component whose repair completes at stage $t$ is not exposed to $\xi^t$); observing $a^t$, the operator dispatches and chooses which failed components to repair, $r^t$; and the repairs, together with $\xi^{t+1}$, determine $a^{t+1}$:
\[
  \xi^t \;\longrightarrow\; a^t \;\longrightarrow\;
  (\text{dispatch } d^t,\ \text{repair } r^t),
  \qquad
  (a^t,\ r^t,\ \xi^{t+1}) \;\longrightarrow\; a^{t+1}.
\]

\paragraph{State variables.} The \emph{component availability vector} $a^t = (a_g^t, a_e^t) \in \{0,1\}^{N_G+N_E}$ records, for each component $k$, whether it is operational ($a^t[k]=1$) or failed ($a^t[k]=0$) at stage $t$, accounting for any failure events and repair completions that have occurred up to and including stage $t$. As a modeling device, all components are operational and no repair is pending before the horizon, $a^0 = \mathbf{1}$ and $r^0 = \mathbf{0}$; the availability $a^1$ then follows from the first failure draw, so components may already be failed at $t=1$.

\paragraph{Stochastic dynamics.} At each stage $t \in \{1,\ldots,T\}$, the failure event vector $\xi^t \in \{0,1\}^{N_G+N_E}$ is revealed. For every component, $\xi^t[k]\sim\operatorname{Bernoulli}(\lambda^t[k])$ independently across time and components, regardless of its state; by~\eqref{eq:transition}, a shock has an effect only on a component that is operational at $t-1$, so $\lambda^t[k]$ is the failure probability of an operational component. The failure process $\xi^1,\ldots,\xi^T$ is the only source of randomness in the model: the load profile $p_l^t$ is exogenous, deterministic, and known to the operator at every stage. The availability at stage $t$ is then updated component-wise according to:
\begin{equation}
  a^t[k] \;=\;
  \begin{cases}
    0 & \text{if } a^{t-1}[k]=1 \text{ and } \xi^t[k]=1
        \quad\text{(new failure)}, \\[4pt]
    1 & \text{if } a^{t-1}[k]=0,\;\; t \geq \tau+1,\;\;\text{and}\quad
        r^{t-\tau}[k] = 1
        \quad\text{(repair after delay }\tau\text{)}, \\[4pt]
    a^{t-1}[k] & \text{otherwise,}
  \end{cases}
  \label{eq:transition}
\end{equation}
where $r^s[k] := 0$ for all $s \leq 0$ (no repairs before the planning horizon). The condition $r^{t-\tau}[k] = 1$ means that the repair was ordered $\tau$ stages earlier. For $\tau=1$ this reduces to $r^{t-1}[k]=1$, which is the case used in the experiments of Section~\ref{sec:experiments}; for general $\tau$, the dynamic-programming state must additionally track repairs in progress through a $\tau$-position pipeline; the MIP formulations used for SAA, SDDiP, and the training labels include this pipeline, whereas the inputs of our learned policy cover $\tau=1$ only.

\paragraph{Stage decisions.} Given $a^t$, the decisions at stage $t$ are:
\begin{itemize}[nosep]
\item \emph{Dispatch}: $(p_g^t, \theta^t, f^t, d^t)$ solving the DC-OPF~\eqref{eq:dc_opf} with availability $(a_g^t, a_e^t)$;
\item \emph{Repair}: $r^t \in \{0,1\}^{N_G+N_E}$ with $r^t[k] \leq 1 - a^t[k]$ (only failed components can be repaired).
\end{itemize}

\paragraph{Stage cost.} The stage cost sums dispatch cost and repair fees:
\begin{equation}
  c_t(p_g^t, d^t, r^t)
  = c^\top p_g^t + \alpha\,\mathbf{1}^\top(p_l^t - d^t)
  + \beta^\top r^t,
  \label{eq:stage_cost}
\end{equation}
where $\beta = (\beta_g, \beta_e)$ collects the component repair costs, $\beta_g \in \R^{N_G}_{\ge0}$ and $\beta_e \in \R^{N_E}_{>0}$ (only a generator with zero capacity has zero repair cost).

\paragraph{The multistage problem.} The repair-scheduling problem is the multistage stochastic mixed-integer program (MSMIP)
\begin{equation}
  \min_\pi \;
  \E\!\left[\sum_{t=1}^T \bigl(c^\top p_g^t + \alpha\,\mathbf{1}^\top(p_l^t - d^t)
    + \beta^\top r^t\bigr)\right]
  \label{eq:msmip}
\end{equation}
subject to, at every stage $t$, the constraints of~\eqref{eq:dc_opf} at availability $a^t$, the dynamics~\eqref{eq:transition} with $a^0=\mathbf{1}$ and $r^0=\mathbf{0}$, and binary repairs $r^t \in \{0,1\}^{N_G+N_E}$ with $r^t \leq \mathbf{1}-a^t$. The policy $\pi$ is non-anticipative: each $d^t$ and $r^t$ is $\sigma(\xi^1,\ldots,\xi^t)$-measurable. The budgeted variant below adds~\eqref{eq:budget}, and Section~\ref{ec:model} states both programs in full.

\paragraph{Extension: a per-stage repair budget.} In practice, the money available for repairs at each stage is often limited, so we also consider a variant of~\eqref{eq:msmip} with a per-stage repair budget. Let $B^t > 0$ be an exogenously given repair budget at stage $t$ (\emph{not} accumulated or carried over); the model then gains the single coupling constraint
\begin{equation}
  \beta^\top r^t \;\le\; B^t, \qquad \forall\,t.
  \label{eq:budget}
\end{equation}
Because the budget is not carried over, it introduces no new state variable, and the state $(a^t,r^t)$ carried to the next stage (Section~\ref{subsec:decoupling}) remains unchanged; however, it \emph{couples} the per-stage repair decisions, which must now be jointly affordable. This single change is sufficient to break the component-wise separability exploited in Section~\ref{subsec:decoupling}: we treat the budgeted variant as an extension of that section and use it to illustrate the general, non-separable case of the framework.

\paragraph{Failure probability model.} The time-varying failure probability $\lambda^t[k]$ of component $k$ is specified by a synthetic log-linear logit model that emulates the non-constant temporal risk profile of a severe-weather event:
\begin{equation}
  \logit(\lambda^t[k]) = A[k]\ln(t+K) - B[k]\,t,
  \label{eq:logit}
\end{equation}
with $K > 0$ a global shift constant. The parameters $A[k]$ and $B[k]$ are set so that $\lambda^t[k]$ peaks at stage $t^*[k]$ with value $p_\text{max}[k]$; Section~\ref{ec:model} gives the expressions, the condition on $K$, and the resulting profiles. With $t^*=1$ and $K=1.5$ for every component, as in our experiments, each failure probability is above $1/2$ up to stage 4 and below it from stage 5 on, whatever its peak value. This specification is a synthetic modeling choice and is not calibrated to meteorological data. The DFL framework and its training procedure do not depend on it: the failure probabilities enter only through the sampled scenarios.

\subsection{Structural analysis: decomposition choice and the threshold policy}
\label{subsec:decoupling}

\begin{proposition}[Dispatch--repair decoupling]
\label{prop:decoupling} Under model~\eqref{eq:msmip}, the dispatch decisions $(p_g^t, \theta^t, f^t, d^t)$ do not appear in the state transition~\eqref{eq:transition}. Consequently, the next-stage availability $a^{t+1}$ and the continuation value depend on stage-$t$ decisions only through the repair vector $r^t$. Moreover, for fixed $a^t$, the optimal dispatch cost $Q_t^0(a^t)$ is independent of $r^t$.
\end{proposition}

\begin{proof}
The dynamics~\eqref{eq:transition} involve only availabilities, shocks, and past repairs, so the dispatch variables do not enter them and $a^{t+1}$ depends on stage-$t$ decisions only through $r^t$. For fixed $a^t$, the DC-OPF~\eqref{eq:dc_opf} contains no repair variable, so its optimal value $Q_t^0(a^t)$ does not depend on $r^t$.
\end{proof}

Proposition~\ref{prop:decoupling} makes $u_t=r^t$ the minimal natural decomposition: the transition depends only on $r^t$, while the stage cost is linear in it, so Assumptions~\ref{assump:state_dep} and~\ref{assump:linear_cost} hold with $\delta_t^0=(p_g^t,\theta^t,f^t,d^t)$. Larger decompositions containing dispatch variables remain valid but unnecessarily embed the DC-OPF in the learning layer. Since $Q_t^0$ does not depend on $r^t$, the surrogate~\eqref{eq:surrogate_u} reduces to a linear problem in only $N_G+N_E$ repair variables, making the FY loss directly applicable (Remark~\ref{rem:fy_connection}); component-wise separability then yields the threshold policy below.

\paragraph{Value convention.} Section~\ref{subsec:msp} defines $\mathcal{V}_{t+1}$ before the shock of stage $t+1$. In the repair problem with $\tau=1$, the state carried from stage $t$ to stage $t+1$ is $x_{t+1}=(a^t,r^t)$: it is observed before $\xi^{t+1}$ and depends on the stage-$t$ decisions only through $r^t$, as Assumption~\ref{assump:state_dep} requires. We write $W_t(a)$ for the optimal expected cost from stage $t$ onwards given the post-shock availability $a^t=a$. With $\mathcal{R}_t(a) := \{r \in \{0,1\}^{N_G+N_E} : r \le \mathbf{1}-a\}$, intersected with~\eqref{eq:budget} when the budget applies, and $a' := a \odot (\mathbf{1}-\xi^{t+1}) + r$ the next availability under~\eqref{eq:transition},
\begin{equation}
\begin{aligned}
  \bar{\mathcal{V}}_t(a, r) &:= \E\bigl[\,W_{t+1}(a') \mid a, r\,\bigr],\\
  W_t(a) &= Q_t^0(a) + \min_{r \in \mathcal{R}_t(a)} \bigl\{\beta^\top r + \bar{\mathcal{V}}_t(a, r)\bigr\},
  \qquad
  W_{T+1} \equiv 0.
\end{aligned}
  \label{eq:exp_cont_value}
\end{equation}
Thus $\bar{\mathcal{V}}_t(a,r)$ is the pre-shock value $\mathcal{V}_{t+1}(x_{t+1})$ of Section~\ref{subsec:msp} at $x_{t+1}=(a,r)$, and the optimal expected cost is $\mathcal{V}_1\bigl((\mathbf{1},\mathbf{0})\bigr) = \E[W_1(\mathbf{1}-\xi^1)]$. It is this \emph{expected} continuation value that we approximate linearly (the specialization of~\eqref{eq:linear_approx} to the repair problem):
\begin{equation}
  \bar{\mathcal{V}}_t(a^t, r^t) \approx h_t(a^t)
    - (\bmu^t)^\top r^t,
  \label{eq:rv_approx}
\end{equation}
where $\bmu^t = (\bmu_g^t, \bmu_e^t) \in \R^{N_G+N_E}$ are the
continuation-value coefficients of the repairs: in~\eqref{eq:rv_approx}, $\bmu^t[k]$ plays the role of the expected future cost reduction from initiating a repair of component $k$ at stage $t$.

In the repair problem, the stage surrogate~\eqref{eq:surrogate_u} becomes
\begin{equation}
  \min_{r^t}\; (\beta - \bmu^t)^\top r^t
  \quad\text{s.t. } r^t[k] \leq 1 - a^t[k],\;\; r^t[k] \in \{0,1\}.
  \label{eq:surrogate_repair}
\end{equation}
Its objective is $\sum_k (\beta[k]-\bmu^t[k])\,r^t[k]$ and each constraint involves a single component, so the problem separates across components. Its optimal solution, with ties broken toward not repairing, is the \emph{threshold repair policy}: repair component $k$ if and only if its estimated future benefit exceeds its immediate cost and it is currently failed, that is
\begin{equation}
  r^{t*}[k] =
  \begin{cases}
    1 & \text{if } \bmu^t[k] > \beta[k] \text{ and } a^t[k] = 0, \\
    0 & \text{otherwise}
  \end{cases}
  \label{eq:threshold}
\end{equation}
(Section~\ref{ec:proofs} gives the formal argument).

\paragraph{Neural network predictor.} A feedforward network $\pi_{\bpsi} : \R^{N_G+N_E+1} \to \R^{N_G+N_E}$ maps the feature vector $(a_g^t, a_e^t, t/T)$ (availability plus normalized stage index) to $\bmu^t$. For convenience, the implementation outputs the effective scores $\eta^t = \bmu^t - \beta$ of Remark~\ref{rem:fy_connection} directly, which is equivalent since $\beta$ is fixed, and~\eqref{eq:threshold} repairs a failed component $k$ if and only if $\eta^t[k] > 0$, with no further optimization.\\
Because the surrogate~\eqref{eq:surrogate_repair} is separable, the smoothed decision map~\eqref{eq:perturbed} is available in closed form. Our unbudgeted runs perturb with standard logistic noise, for which
\begin{equation}
  \hat{r}_\varepsilon(\eta^t)[k]
  = \sigma\!\left(\frac{\eta^t[k]}{\varepsilon}\right)
    \cdot \mathbf{1}[a^t[k] = 0],
  \label{eq:smooth_threshold}
\end{equation}
where $\sigma(z) = (1+e^{-z})^{-1}$ is the logistic function. The training gradient is $[\hat{r}_\varepsilon(\eta^t) - r^{t*}]$, requiring no solver calls at training time; Section~\ref{ec:impl} compares this choice with Gaussian noise.

\paragraph{Extension: when a budget couples the components.} The budgeted variant~\eqref{eq:budget} is exactly a non-separable instance. With the per-stage constraint $\beta^\top r^t \le B^t$, the separable surrogate~\eqref{eq:surrogate_repair} becomes a $0/1$ knapsack over the failed components,
\begin{equation}
  r^{t*}\in\arg\max_{r}\;(\bmu^t-\beta)^\top r
  \quad\text{s.t.}\quad \beta^\top r\le B^t,\;\; r[k]\le 1-a^t[k],\;\; r[k]\in\{0,1\},
  \label{eq:surr_knap}
\end{equation}
with item value $(\bmu^t[k]-\beta[k])$, weight $\beta[k]$, and capacity $B^t$ (at deployment the values are nominal, so negative-value items, $\bmu^t[k]<\beta[k]$, never enter; in training they are perturbed first, so any failed component may enter). Its feasible set
\begin{equation}
  \mathcal Y_t=\bigl\{r\in\{0,1\}^{N_G+N_E}:\ r\preceq \mathbf 1-a^t,\ \ \beta^\top r\le B^t\bigr\}
  \label{eq:Yfeas_budget}
\end{equation}
is no longer a box, so the perturbed optimum no longer factors across components and the threshold rule~\eqref{eq:threshold} and closed-form smoothing~\eqref{eq:smooth_threshold} no longer hold. The Fenchel--Young loss~\eqref{eq:fy} keeps its form with effective score $\eta^t:=\bmu^t-\beta$, and the gradient is still $\hat r_\varepsilon(\eta^t)-r^{t*}$; but the smoothed decision map is now the perturbed knapsack optimum
\begin{equation}
  \hat r_\varepsilon(\eta^t)=\E_{Z\sim\mathcal N(0,I)}\!\Bigl[\arg\max_{r\in\mathcal Y_t}(\eta^t+\varepsilon Z)^\top r\Bigr],
  \label{eq:perturbed_knap}
\end{equation}
which has no componentwise closed form (because $\mathcal Y_t$ is not a box) and is estimated by Monte Carlo perturbation: one knapsack solve per sample, exactly the general estimator of Section~\ref{subsec:fy}, with only the objective coefficients perturbed while the weights $\beta$ and the capacity $B^t$ stay fixed. The single-scenario training labels and all baselines (Section~\ref{sec:experiments}) also respect~\eqref{eq:budget}. A slack budget returns the deployment decision to the threshold rule~\eqref{eq:threshold} (Section~\ref{ec:proofs} gives the exact conditions, including those for the smoothed map), so a binding budget is what makes the extension non-trivial: it removes the separability specific to the power-network problem and requires the general estimator of Section~\ref{subsec:fy}.

\section{Computational Experiments}
\label{sec:experiments}

\subsection{Experimental Setup}
\label{subsec:setup}

\paragraph{Test networks.} Two benchmark networks are used: a synthetic 5-bus network ($N_G=2$, $N_E=5$), which serves as a small validation instance, and the IEEE 14-bus network, case14 ($N_G=5$, $N_E=20$). All experiments run over a $T=14$ stage horizon with repair delay $\tau=1$. In the DC-OPF~\eqref{eq:dc_opf}, the voltage-angle bound is set to $\theta^{\max}=\pi$ in all experiments, so the big-$M$ constants are $M_j = 2\pi |b_j|$; the generation lower limits are $\underline{p}_g = \mathbf{0}$ (no must-run generation). Demand is constant over the horizon, $p_l^t = p_l$ at every stage, equal to the nominal bus loads of the network data, and each stage is costed as one hour of operation, so a shortfall in MW corresponds to the same amount in MWh of unserved energy.

\paragraph{Failure parameters.} Components are assigned to two vulnerability classes with equal probability (50\% each): \emph{standard} ($p_\text{max} \sim \mathcal{U}[0.55, 0.70]$)
and \emph{elevated} ($p_\text{max} \sim \mathcal{U}[0.80, 0.90]$), with peak failure stage $t^*=1$ and global shift $K=1.5$. The value of lost load is $\alpha=\$1{,}000$/MWh. Baseline repair costs are capacity-proportional: generator mean $\$550$ and line mean $\$350$. Because these failures are synthetic and independent across components, and their probabilities peak already at the first stage, the experiments below are algorithmic stress tests that compare solution methods under one deliberately severe failure regime.

\paragraph{DFL training.} The neural network has three hidden layers and is trained with Adam and early stopping at perturbation temperature $\varepsilon=0.1$; Section~\ref{ec:impl} lists all settings and documents the stop-gradient form of the loss used in the implementation. Training labels come from single-scenario MIPs solved to a relative gap of $10^{-6}$. In the budgeted variant the perturbed knapsack oracle~\eqref{eq:perturbed_knap} is estimated with $M=5$ Monte Carlo samples per training example. We sweep 15 training sizes $n$ from 25 to 120{,}000 single-scenario MIPs. Each MIP contributes $T=14$ stage-level examples, and the training sets are nested: each size uses the first $n$ scenarios of one common sample. The stage-level examples are split at random into 80\%~training and 20\%~validation, and the validation loss drives the learning-rate schedule, early stopping, and the choice of the restored checkpoint. Every size is trained from five weight-initialization seeds, which change only the initial weights; the data, the split, the dropout masks, and the batch order are shared.
\paragraph{Baselines.} Five baselines are evaluated alongside DFL.
\emph{NoRepair} never initiates a repair, and \emph{AlwaysRepair} repairs every failed component at every stage except the last, where a repair could not complete within the horizon; these are the two trivial heuristics.
\emph{SAA} solves a single extensive-form MILP on a sampled binary scenario tree over the 14 stages and is evaluated by rolling its node policy out on the common evaluation paths with nearest-child matching. \emph{SDDiP} runs Lagrangian cut generation through \texttt{SDDP.jl}~\citep{dowson2021sddpjl} until the shared wall-clock limit, and its final policy is evaluated on 500 simulations drawn from its own stagewise distribution rather than on the common evaluation paths. Section~\ref{ec:comparators} gives both specifications.\\
Because the two trivial heuristics are near-optimal only in the extreme cost regimes (Section~\ref{subsec:sens}), we add a non-anticipative baseline that plans ahead: \emph{Near-sighted MPC} (Model Predictive Control), a certainty-equivalent, finite-horizon model-predictive controller. At each stage it rolls the failure process forward at its mean over a window of $W=3$ stages, scoring the resulting fractional availabilities with a DC-OPF whose line and generator parameters are scaled by them, optimizes the repair schedule over that window, executes only the current-stage repairs, and replans after observing the next failures. We call it ``near-sighted'' since it assigns zero value to anything beyond its planning window (no terminal cost-to-go or continuation value). Its schedule search is exhaustive for up to seven simultaneously failed components and samples schedules beyond that. It is non-anticipative by construction and uses no learning; Section~\ref{appendix:mpc} presents the algorithm in full and details how it differs from a complete stochastic MPC.\\
Because the 5-bus instance has only seven failure-prone components, we additionally compute its exact optimal expected value by finite-horizon dynamic programming and use it as a benchmark. For the empirical policy comparison, all methods except SDDiP, which is evaluated internally as described above, are tested on the same 500 \emph{evaluation paths}, sampled with a seed different from that of the training scenarios. These paths are separate from model fitting, but we also used them to choose the operating scales (Section~\ref{subsec:sens}) and, for each configuration, the reported training size (the one with the lowest mean cost).

\paragraph{Compute environment.} Every run, for every method, used the same type of node on a high-performance computing cluster (AMD EPYC 7H12 CPUs; 64 cores and 112\,GB of memory requested per job; no GPUs) under a wall-clock limit of 117 hours (roughly 4.9 days) per job; Section~\ref{ec:compute} describes how each method used these resources. The code, with its solver settings, job scripts, and seeds, is publicly available.

\subsection{Cost-Regime Analysis and Parameter Selection}
\label{subsec:sens}

Depending on the ratio between the load-shedding penalty $\alpha$ and the component repair costs $\beta$, the problem tends toward one of two trivial extremes. When repairs are nearly free, repairing every failed component immediately is near-optimal (\emph{AlwaysRepair}); when repairs are prohibitively expensive, never repairing dominates (\emph{NoRepair}). The informative regime lies in between, where the operator must decide \emph{which} components to repair and \emph{when}. To locate it, we fix $\alpha=\$1{,}000$/MWh and scale all repair costs by a common factor $s$.

\begin{figure}[tb]
{\centering
\begin{minipage}[b]{0.48\linewidth}\centering
  \includegraphics[width=\linewidth,height=4.5cm,keepaspectratio]{%
    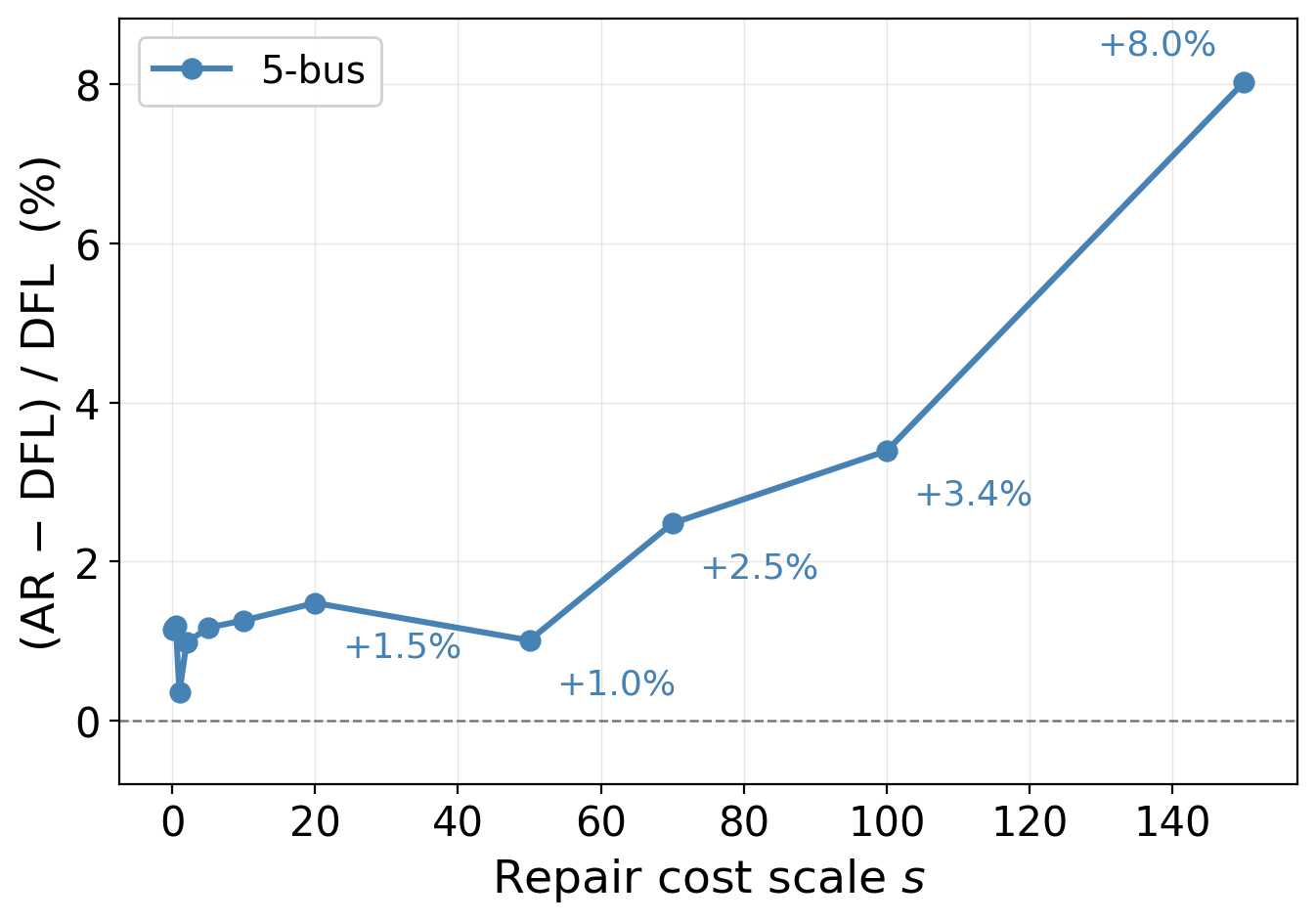}
  \\[-3pt]\small\textit{(a) 5-bus ($s=100$ selected)}
 \end{minipage}\hfill
 \begin{minipage}[b]{0.48\linewidth}\centering
  \includegraphics[width=\linewidth,height=4.5cm,keepaspectratio]{%
    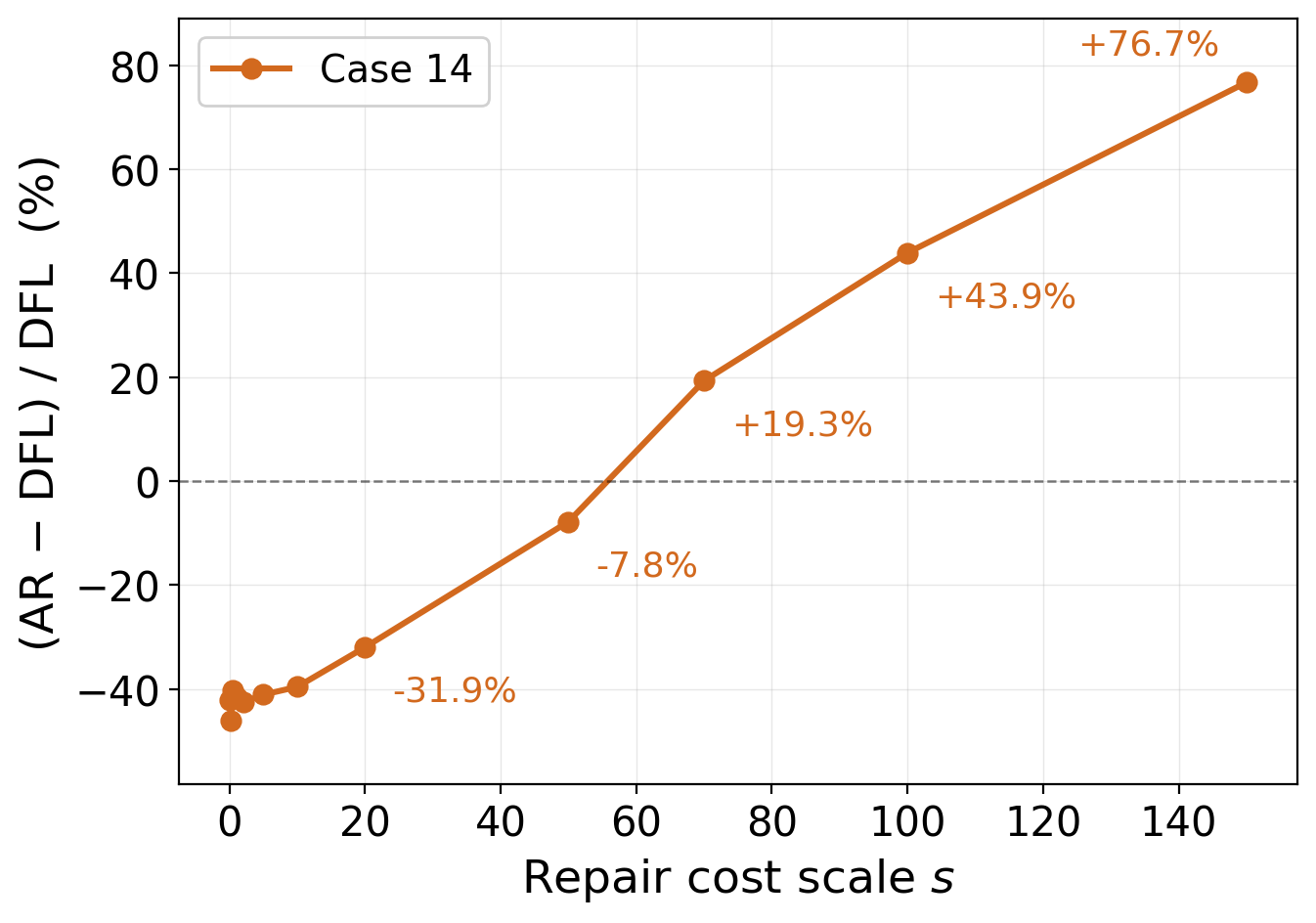}
  \\[-3pt]\small\textit{(b) case14}
 \end{minipage}\par}
\caption{Relative advantage of DFL over AlwaysRepair vs.\ repair cost scale $s$, unbudgeted problem.\label{fig:cost_vs_scale}}
{\small\textit{Note.} Vertical axis: $(\bar c_\text{AR}-\bar c_\text{DFL})/\bar c_\text{DFL}\times100\%$; positive means DFL is cheaper. Each point is a single-seed DFL model trained on $n=120{,}000$ scenarios and evaluated on the 500 evaluation paths, so the values differ from the five-seed means of Table~\ref{tab:dfl}. The case14 operating scales $s=45$ and $s=55$ lie between the sweep points $s=20$, $50$, and $70$.\par}

\end{figure}

Figure~\ref{fig:cost_vs_scale} sweeps $s\in[0.1,150]$ on the unbudgeted problem and reports the relative advantage of DFL over AlwaysRepair. On 5-bus, DFL is ahead at every scale tested, with a margin that stays small while repairs are cheap and grows once they become expensive. On case14, AlwaysRepair is strictly better over most of the range; the advantage changes sign between the sweep points $s=50$ and $s=70$, near $s\approx56$ by linear interpolation, and widens quickly beyond them. On both networks, DFL's advantage over AlwaysRepair thus grows with the repair-cost scale.

We choose the operating scales from the \emph{unbudgeted} problem. A binding per-stage cap~\eqref{eq:budget} removes feasible repair schedules and forces the operator to choose among failed components, so a scale at which repair selection already matters without the cap remains relevant once the cap is imposed. We fix
\begin{equation}
  s_\text{5-bus}=100,\qquad s_\text{case14}\in\{45,\,55\},
  \label{eq:opscales}
\end{equation}
and all experiments in Sections~\ref{subsec:sample_complexity}--\ref{subsec:time} use these scales, with every method evaluated at the same per-network scale. For 5-bus, $s=100$ lies in the range where DFL's advantage over AlwaysRepair grows with $s$. For case14, we choose two scales around the unbudgeted crossover: $s=45$, where repairs are cheaper relative to the load-shedding penalty and selective repair has less room to pay off, and $s=55$, close to the crossover. Section~\ref{subsec:dfl_vs} compares the unbudgeted results at both scales, while Section~\ref{subsec:budget_results} reports how the comparison at $s=45$ changes once the cap restricts the repairs that can be made at each stage.

\subsection{Empirical Learning Curves}
\label{subsec:sample_complexity}

We now examine how DFL's cost on the evaluation paths changes with the size of the training set at the operating scales in~\eqref{eq:opscales}. For each network we train DFL at 15 sizes from 25 to 120{,}000 single-scenario MIPs and evaluate each trained model on the 500 evaluation paths; the budgeted $s=55$ sweep is capped at $n=80{,}000$. To distinguish the effect of training-set size from variation due to weight initialization, we train five independently initialized models at each size.

\begin{figure}[tb]
{\centering
\begin{minipage}{\linewidth}
 \begin{minipage}[b]{0.48\linewidth}\centering
  \includegraphics[width=\linewidth,height=4.6cm,keepaspectratio]{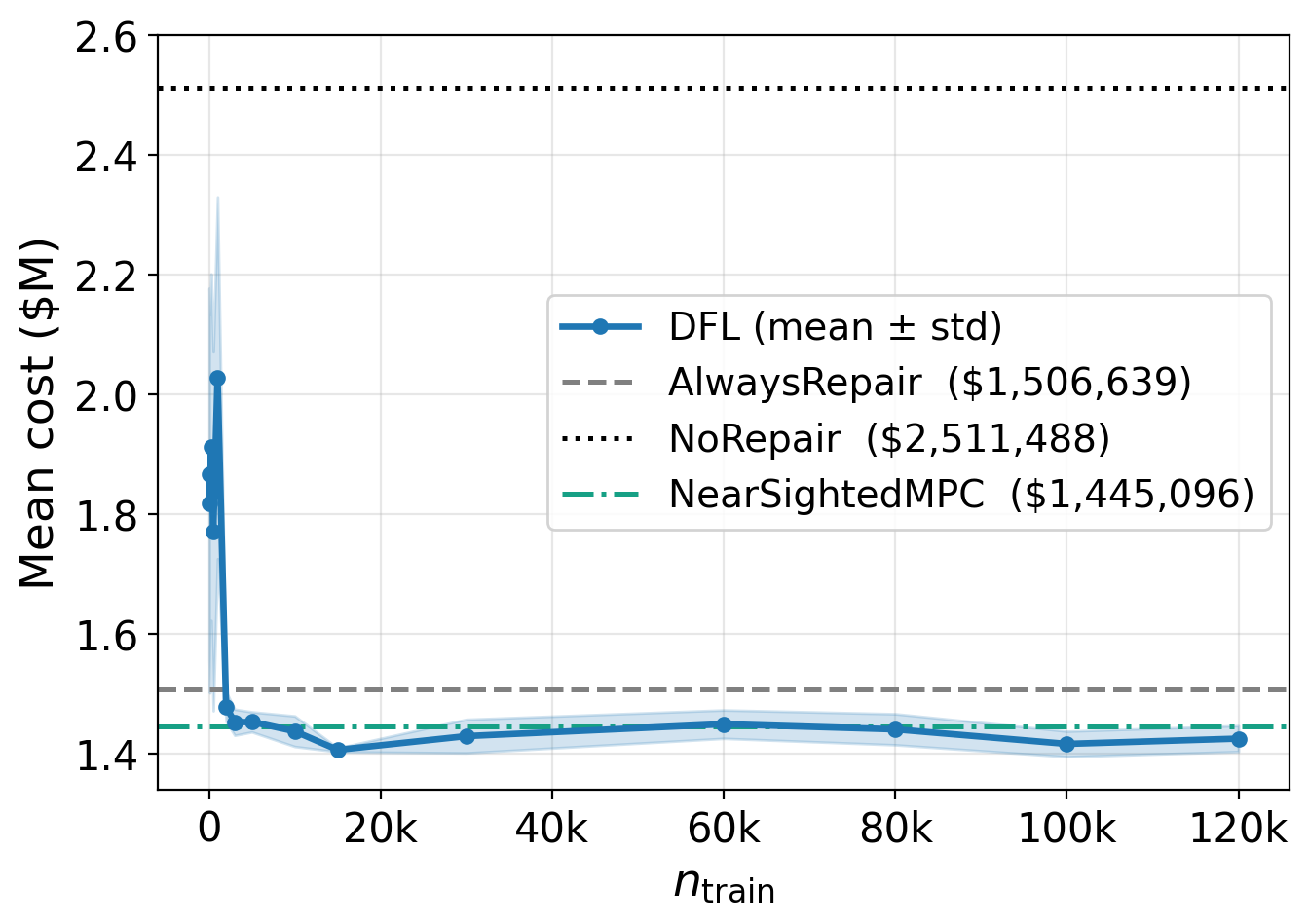}
  \\[-3pt]\small\textit{(a) 5-bus ($s=100$), unbudgeted}
 \end{minipage}\hfill
 \begin{minipage}[b]{0.48\linewidth}\centering
  \includegraphics[width=\linewidth,height=4.6cm,keepaspectratio]{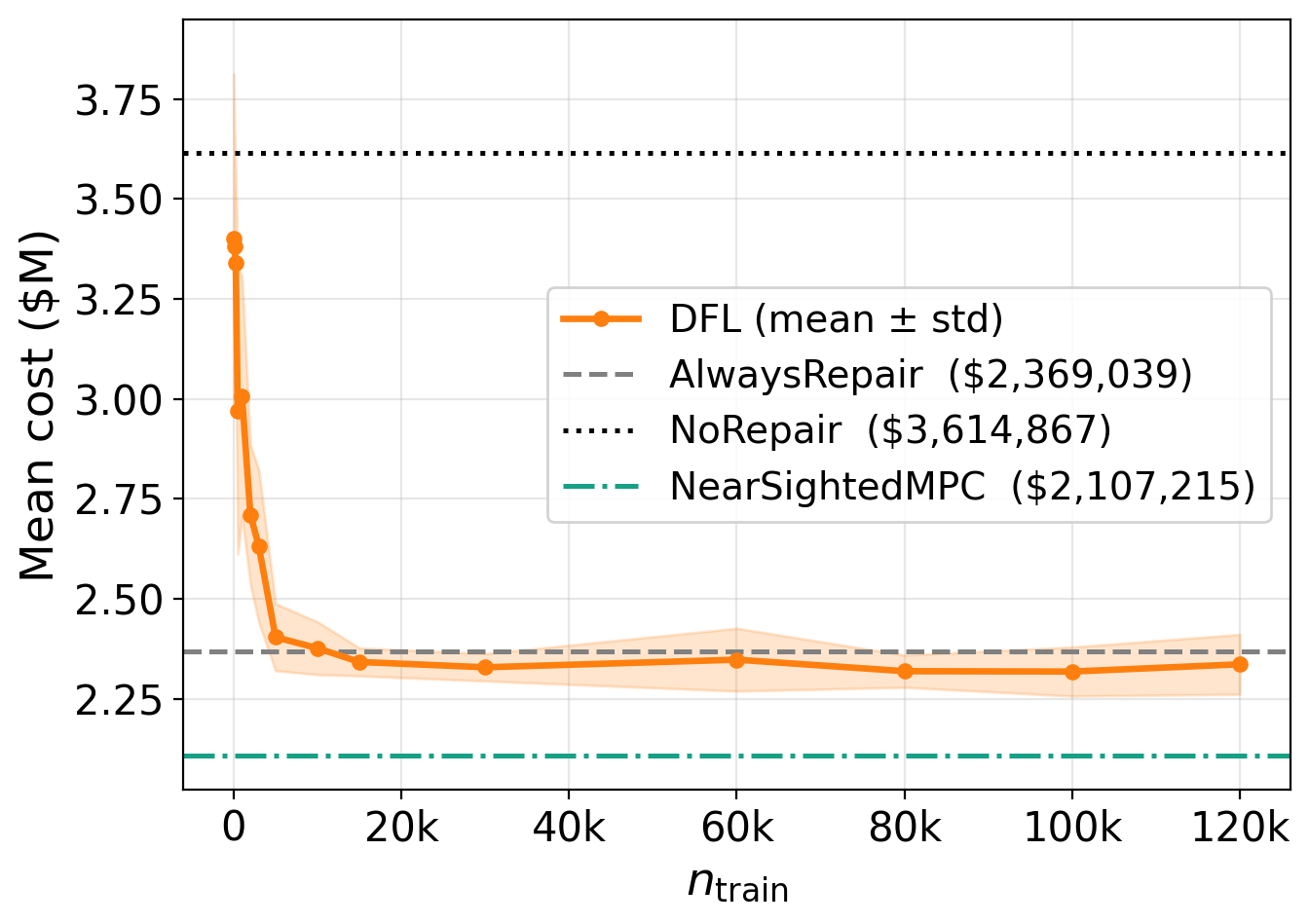}
  \\[-3pt]\small\textit{(b) case14 ($s=55$), unbudgeted}
 \end{minipage}\\[6pt]
 \begin{minipage}[b]{0.48\linewidth}\centering
  \includegraphics[width=\linewidth,height=4.6cm,keepaspectratio]{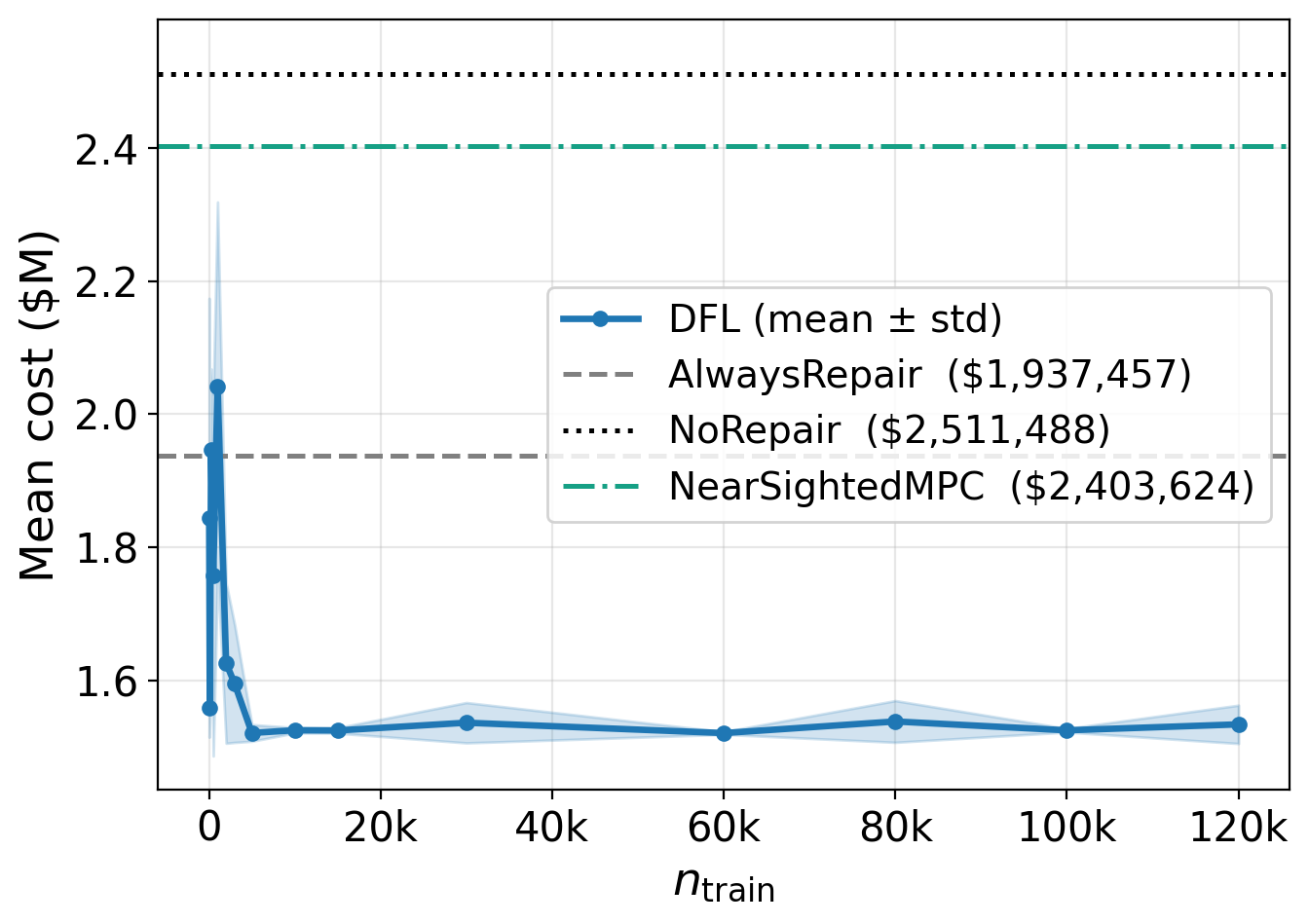}
  \\[-3pt]\small\textit{(c) 5-bus ($s=100$), budgeted}
 \end{minipage}\hfill
 \begin{minipage}[b]{0.48\linewidth}\centering
  \includegraphics[width=\linewidth,height=4.6cm,keepaspectratio]{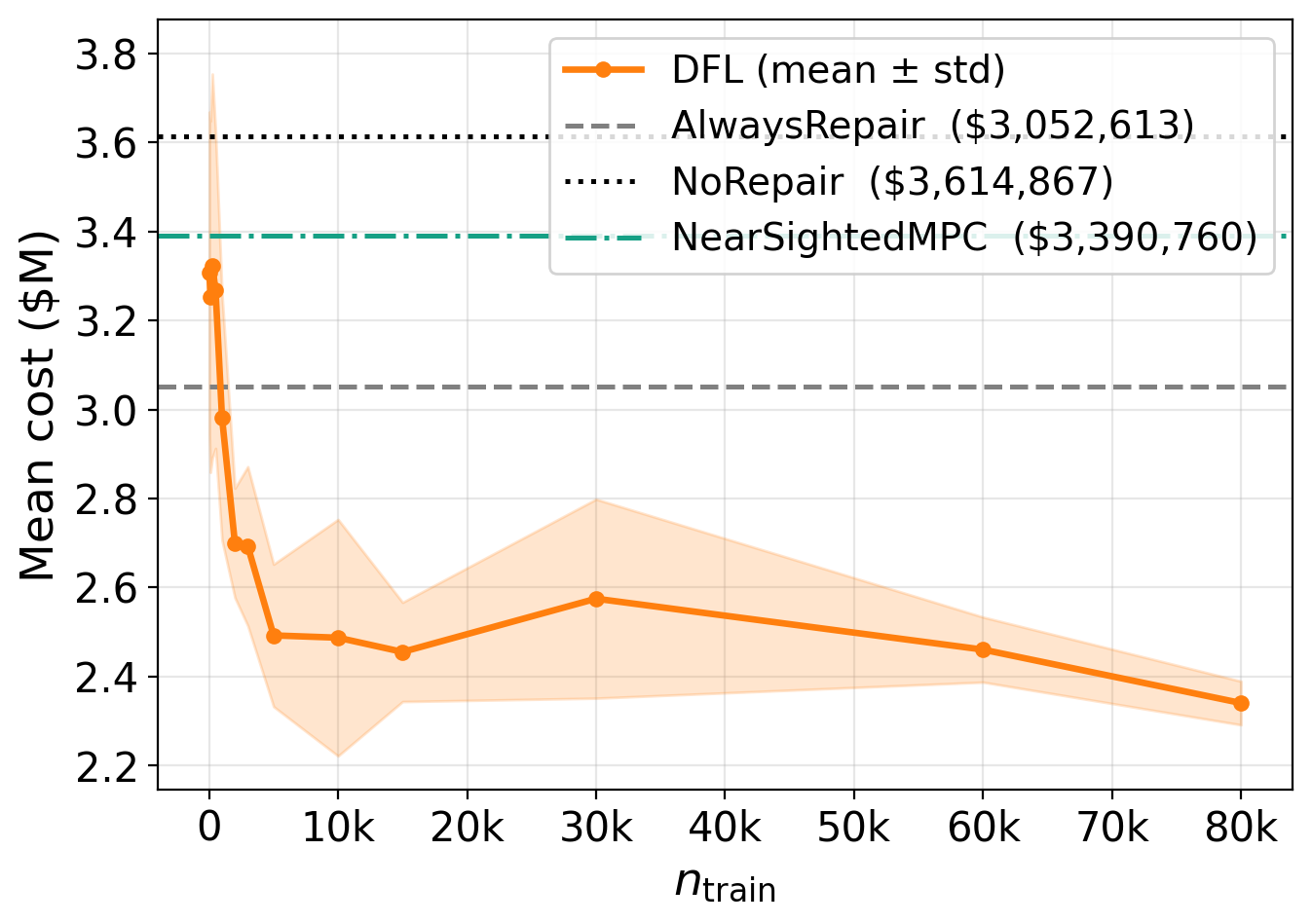}
  \\[-3pt]\small\textit{(d) case14 ($s=55$), budgeted}
 \end{minipage}
 \end{minipage}\par}
\caption{DFL mean cost on the evaluation paths ($\pm$std across five seeds) vs.\ training size.\label{fig:dfl_convergence}}
{\small\textit{Note.} Panels (c) and (d) use the per-stage budget ($\varphi=0.3$, Section~\ref{subsec:budget_results}) and are trained through the perturbed-knapsack oracle; the budgeted $s{=}55$ sweep is capped at $n=80{,}000$. Dashed: AlwaysRepair; dotted: NoRepair; dash-dotted: Near-sighted MPC.\par}

\end{figure}

Figure~\ref{fig:dfl_convergence} shows the DFL mean cost on the evaluation paths $\pm$ one standard deviation across the five seeds, as a function of $n_\text{train}$, without the budget (panels a and b) and with it (panels c and d).

Most of the improvement occurs at relatively small training sizes. On \emph{5-bus} ($s=100$) for instance, performance varies considerably below $n\approx1{,}000$, but DFL moves below AlwaysRepair at $n=2{,}000$ and changes little beyond $n\approx15{,}000$; over this range, the standard deviation across seeds remains below $\$30$k. Section~\ref{subsec:dfl_vs} compares this policy with the baselines in greater detail.

\subsection{DFL versus Baselines}
\label{subsec:dfl_vs}

Table~\ref{tab:dfl} reports, for each network at its operating scale, the mean cost of every method over the 500 evaluation paths, without the budget (upper panel) and with it (lower panel, Section~\ref{subsec:budget_results}).

\begin{table}[tb]
\caption{Mean cost on the evaluation paths by method at the chosen operating scales, without and with a binding per-stage repair budget ($\varphi=0.3$).\label{tab:dfl}}
{\centering
\small\begin{tabular}{@{}lrrrrrrr@{}}
\toprule
Network & $n_\text{train}$ & NoRepair & AlwaysRepair & MPC & SAA & DFL & Savings \\
\midrule
\multicolumn{8}{@{}l}{\textit{Without budget}} \\
5-bus ($s{=}100$)   & 15{,}000  & 2{,}511{,}488 & 1{,}506{,}639 & 1{,}445{,}096 & 2{,}028{,}014 & \textbf{1{,}406{,}419} & $+6.7\%$ \\
case14 ($s{=}45$) & 30{,}000  & 3{,}614{,}867 & 2{,}098{,}438 & \textbf{1{,}949{,}190} & 3{,}193{,}445 & 2{,}271{,}956 & $-8.3\%$ \\
case14 ($s{=}55$) & 100{,}000 & 3{,}614{,}867 & 2{,}369{,}039 & \textbf{2{,}107{,}215} & 3{,}255{,}115 & 2{,}318{,}269 & $+2.1\%$ \\
\midrule
\multicolumn{8}{@{}l}{\textit{With per-stage budget ($\varphi=0.3$)}} \\
5-bus ($s{=}100$)   & 60{,}000 & 2{,}511{,}488 & 1{,}937{,}457 & 2{,}403{,}624 & 2{,}069{,}534 & \textbf{1{,}520{,}472} & $+21.5\%$ \\
case14 ($s{=}45$) & 15{,}000 & 3{,}614{,}867 & 2{,}842{,}476 & 2{,}600{,}141 & 3{,}240{,}817 & \textbf{2{,}235{,}973} & $+21.3\%$ \\
case14 ($s{=}55$) & 80{,}000 & 3{,}614{,}867 & 3{,}052{,}613 & 3{,}390{,}760 & 3{,}214{,}677 & \textbf{2{,}340{,}615} & $+23.3\%$ \\
\bottomrule
\end{tabular}\par}
\smallskip
{\small\textit{Note.} Mean total cost (\$) over the 500 evaluation paths; bold marks the cheapest policy in each row. DFL is the five-seed mean at its best training size; with the budget, its across-seed std is $\$2{,}011$, $\$68{,}890$, and $\$49{,}014$ (top to bottom). Savings $=(\bar c_\text{AR}-\bar c_\text{DFL})/\bar c_\text{AR}\times100\%$, positive when DFL is cheaper than AlwaysRepair (AR). With the budget, AlwaysRepair repairs as many failed components as the budget allows, cheapest first, except at the last stage; the schedule filter of Near-sighted MPC receives the cap; and DFL solves the knapsack~\eqref{eq:surr_knap}.\par}

\end{table}

Across the two case14 rows, DFL's saving over AlwaysRepair changes sign: DFL remains behind the strongest baselines at $s=45$ but moves ahead of AlwaysRepair at $s=55$. On 5-bus, DFL is the cheapest policy evaluated, $6.7\%$ below AlwaysRepair, though Near-sighted MPC is close behind at $2.8\%$ above DFL. The exact optimal expected cost, computed by dynamic programming, is $\$1{,}382{,}625$, so the DFL mean lies $1.7\%$ above it. On case14 at $s=45$, DFL is not cheapest: Near-sighted MPC and AlwaysRepair both undercut it, by $14.2\%$ and $7.6\%$ of the DFL cost. At $s=55$, DFL moves $2.1\%$ below AlwaysRepair, though Near-sighted MPC remains the cheapest policy at that scale; Section~\ref{subsec:budget_results} shows that this ordering reverses once the cap binds.

The cost decomposition shows where the differences arise. On case14 at $s=45$, DFL spends $\$540{,}798$ on repairs, compared with $\$1{,}217{,}706$ for AlwaysRepair, but sheds $1{,}714$\,MWh instead of $857$\,MWh, and the additional load-shedding cost exceeds the saving on repairs. At $s=55$, DFL spends $\$517{,}535$ against $\$1{,}488{,}307$ and sheds $1{,}784$\,MWh against $857$\,MWh. On 5-bus at $s=100$, DFL spends $\$381{,}209$ on repairs, compared with $\$737{,}884$ for AlwaysRepair, while shedding $988$\,MWh instead of $726$\,MWh; here the additional load-shedding cost is smaller than the repair expenditure avoided. In all three unbudgeted configurations, DFL repairs less and sheds more load than AlwaysRepair.\\
Across scenarios, the standard deviation of the DFL cost is lower than that of AlwaysRepair and Near-sighted MPC on 5-bus and on case14 at $s=55$, and comparable at $s=45$.

\subsection{Budgeted Repair Extension: Results}
\label{subsec:budget_results}

We next consider the budgeted variant~\eqref{eq:budget}, in which the per-stage constraint $\beta^\top r^t\le B^t$ couples the repair decisions. As established in Section~\ref{subsec:decoupling}, a binding budget removes the separability of the unbudgeted surrogate. The stage surrogate becomes the $0/1$ knapsack problem~\eqref{eq:surr_knap}, the closed-form rule~\eqref{eq:threshold} is no longer applicable, and training uses the general perturbed-optimizer estimator of Section~\ref{subsec:fy}. We set the stage budget to $\varphi=0.3$ of the cost of repairing all components, so that $B^t=0.3\,\beta^\top\mathbf{1}$ at every stage. The networks, operating scales~\eqref{eq:opscales}, and evaluation paths are those of the unbudgeted experiments, while the labels, the stage surrogate, and the baselines are adapted to the budget. The lower panel of Table~\ref{tab:dfl} reports the results.

The comparison reverses on case14, while DFL's existing advantage on 5-bus becomes larger. On 5-bus, the budget increases the cost of AlwaysRepair much more than that of DFL, raising DFL's saving from $6.7\%$ to $21.5\%$; the DFL mean lies $1.8\%$ above the exact optimal expected cost of $\$1{,}493{,}319$. On case14 at $s=45$, DFL moves from costing $8.3\%$ more than AlwaysRepair in the unbudgeted setting to costing $21.3\%$ less under the budget. At $s=55$, the corresponding saving is $23.3\%$. Near-sighted MPC, whose schedule search must respect the additional constraint, becomes much more expensive under the budget: its cost reaches $\$2{,}403{,}624$ on 5-bus, $\$2{,}600{,}141$ at $s=45$, and $\$3{,}390{,}760$ at $s=55$. DFL reduces these costs by $36.7\%$, $14.0\%$, and $31.0\%$, respectively, and is the least-cost policy in all three budgeted configurations.

The cost decomposition indicates where these gains arise. On case14 at $s=45$, DFL spends $\$379{,}670$ on repairs, compared with $\$945{,}617$ for AlwaysRepair, and sheds slightly less load ($1{,}840$ versus $1{,}878$\,MWh). On 5-bus, it lowers both repair expenditure ($\$342{,}405$ versus $\$553{,}488$) and load shedding ($1{,}144$ versus $1{,}356$\,MWh). At $s=55$, DFL produces nearly the same amount of load shedding ($1{,}889$ versus $1{,}878$\,MWh) while reducing repair expenditure from $\$1{,}155{,}754$ to $\$434{,}735$.

Figure~\ref{fig:dfl_convergence}(c)--(d) presents the budgeted learning curves. Under the budget, DFL has the smallest across-scenario standard deviation of the policies that repair, in all three configurations: \$153{,}482, \$247{,}306, and \$215{,}595 for 5-bus, case14 at $s=45$, and case14 at $s=55$, respectively (averaged over the five seeds), against \$216{,}198, \$298{,}073, and \$306{,}815 for AlwaysRepair.

\subsection{Classical Approaches}
\label{subsec:classical}

We compare DFL with our implementations of SAA and SDDiP in the six configurations defined by the chosen operating scales~\eqref{eq:opscales} and the two budget regimes.

\paragraph{SAA.} On 5-bus, SAA reaches its $10\%$ target gap within hours. In the four case14 configurations, it reaches the wall-clock limit with MIP gaps between $42.5\%$ and $60.1\%$. The sampled-tree incumbent objectives range from $\$1.41$M to $\$2.87$M. When the resulting node policies are rolled out on the 500 evaluation paths using nearest-child matching, their costs rise to $\$2.03$M--$\$3.26$M, above both DFL and AlwaysRepair in all six configurations (Table~\ref{tab:dfl}).

\paragraph{SDDiP.} Within the wall-clock limit, SDDiP completes $782$ and $755$ iterations on 5-bus (unbudgeted and budgeted) and at most $158$ on case14 (Table~\ref{tab:sddip}). Measured against the mean cost of its final policy over 500 simulations from its own stagewise distribution, the final lower bound leaves gaps of $35.3\%$ and $21.5\%$ on 5-bus and between $59.8\%$ and $68.2\%$ on case14, so none of the runs converged within the limit. On 5-bus, exact dynamic programming provides an additional check: the final policy means are $14.6\%$ and $10.6\%$ above the exact optima, while the lower bounds remain $25.9\%$ and $13.2\%$ below them. The lower bound moved slowly throughout: in the unbudgeted 5-bus run it rose from $\$0.40$M after the first iteration to $\$0.73$M after 100 and $\$1.02$M after 782. Because its case14 simulations use their own sampled law, we do not rank SDDiP against the policies evaluated on the common evaluation paths, and we use its gaps only as an approximate measure of convergence.

\subsection{Computational Efficiency}
\label{subsec:time}

At deployment, DFL is also computationally cheap, in every cost regime. A complete 14-stage evaluation path takes DFL about $30$\,ms on 5-bus ($s=100$) and $35$\,ms on case14 ($s=45$), comparable to the trivial heuristics ($25$--$36$\,ms) and less than half the time of the deployed SAA policy ($77$--$79$\,ms). Evaluation requires one network forward pass, the thresholding rule~\eqref{eq:threshold}, and one single-stage LP at each stage, and the budgeted policy also solves one knapsack~\eqref{eq:surr_knap} per stage; no scenario sample is involved, and the cost of a complete path grows linearly with $T$. Near-sighted MPC, in contrast, takes $5.9$\,s and $95.3$\,s per path, $196$ and $2{,}707$ times longer than DFL, because the number of simultaneously failed components enlarges the set of repair schedules it must evaluate. DFL's main cost is incurred offline, during training, which is performed once per configuration; training a single model took at most $17$\,h without the budget and $63$\,h with it, within the 117-hour job limit. Generating the training labels is a further offline cost: without the budget, solving the 120{,}000 case14 single-scenario MIPs took about $19$\,h of wall-clock time with up to 64 parallel solves. The single-scenario MIPs are independent of one another, so label generation can run in parallel, and the whole offline cost is incurred once per configuration: the trained policy is then reused on every new failure realization without further training.

\section{Discussion and Conclusion}
\label{sec:conclusion}

This paper develops a decision-focused learning framework for multistage stochastic programs. Under the conditions of Section~\ref{sec:dfl_general}, the effect of a transition-relevant sub-decision on the future is represented through learned continuation-value coefficients, and the stage surrogate becomes a linear-objective oracle that is trained with a Fenchel--Young loss through perturbed optimizers, without differentiating through the solver; the resulting policy is optimal when the continuation value is affine and its coefficients are recovered (Proposition~\ref{prop:conditional}). The power-network problem serves as a worked example: the unbudgeted variant yields a separable threshold rule, while the budgeted variant illustrates the more general case of a coupled decision set and a combinatorial oracle.

\paragraph{Computational findings.} On the common evaluation paths, DFL is the least-cost policy evaluated in all three budgeted configurations. Without the budget, it is cheaper than AlwaysRepair on 5-bus and on case14 at $s=55$ but not at $s=45$, and Near-sighted MPC is cheaper than DFL in two of the three configurations. On the 5-bus network, the DFL mean lies $1.7\%$ (unbudgeted) and $1.8\%$ (budgeted) above the exact optimum, and at deployment DFL is about as fast as the trivial heuristics. The classical results concern our implementations within the shared budget: the sampled-tree SAA policy with nearest-child matching is not competitive on the evaluation paths, and SDDiP remains far from convergence.

\paragraph{Limitations and future work.} When the true continuation value is not affine in the transition-relevant decision, the learned policy has no certified optimality-gap bound. Its performance may depend on approximation error, interactions among decisions, perfect-information labels, state-space coverage, model capacity, and numerical training error. The empirical learning curves do not provide a finite-sample guarantee, and the non-convex neural parameterization precludes a general claim of global convergence.

The 500 evaluation paths are separate from model fitting but are also used to examine the operating scales and to select the reported training sizes. The results should therefore be read as comparisons on a fixed evaluation set. Moreover, the experiments cover only two networks under a synthetic, independent failure model.

Further work should study richer continuation-value models, policy-gap and sample-complexity guarantees, the effect of perfect-information labels on the deployed policy (for instance, through a generalization procedure such as the one applied in \citet{Maximilian_DFL_EMS}), non-linear surrogates, correlated failures, and applications beyond power-network repair. A common-scenario evaluation of a retained SDDiP policy would also sharpen the comparison with the classical methods.

\section*{Acknowledgments}
D. Krame Kadurha and A. Zocca are supported by the NWO grant VI.Vidi.233.247, \emph{Power Network Optimization in the Age of Climate Extremes} (\url{https://doi.org/10.61686/GOOEL09973}). We thank SURF (www.surf.nl) for the support in using the National Supercomputer Snellius. We also used generative AI tools in preparing this manuscript: Claude Code (Claude Opus 5.5, Anthropic) assisted with coding and data analysis, and Claude (Opus 5.5, Anthropic) and ChatGPT (GPT-5.6 Luna, OpenAI) assisted with language editing. The authors reviewed and verified all AI-assisted output and take full responsibility for the content of this article.

\bibliographystyle{informs2014}
\bibliography{references}

\clearpage
\appendix
\counterwithin{equation}{section}
\counterwithin{figure}{section}
\counterwithin{table}{section}
\counterwithin{theorem}{section}
\counterwithin{proposition}{section}
\counterwithin{assumption}{section}
\counterwithin{remark}{section}

\section*{Appendix}

\noindent This appendix collects the material referred to in the main paper. Equations, sections, propositions, and theorems numbered without a letter prefix refer to the main paper.

\section{Fenchel--Young Background and Proofs}
\label{ec:fy}

\subsection{The Fenchel--Young loss}

Let $\mathcal{Y} \subseteq \R^m$ be a compact feasible set (which may be binary, discrete, continuous, or mixed), $\mathcal{C} = \conv(\mathcal{Y})$ its convex hull, and $\Omega: \mathcal{C} \to \R$ a proper, lower-semicontinuous, strongly convex regularizer. Define the regularized oracle and its Fenchel conjugate:
\begin{equation}
  \hat{u}_\Omega(\eta) = \arg\max_{u \in \mathcal{C}}
    \{\eta^\top u - \Omega(u)\},
  \qquad
  \Omega^*(\eta) = \max_{u \in \mathcal{C}}\{\eta^\top u - \Omega(u)\}.
  \label{eq:oracle}
\end{equation}
The Fenchel--Young loss associated with $\Omega$ is $L_\Omega(\eta;u) = \Omega^*(\eta) + \Omega(u) - \eta^\top u$, as in~\eqref{eq:fy}.

\begin{proposition}[\citealt{blondel2020learning}]
\label{prop:fy_props} Let $\Omega$ be strongly convex on $\mathcal{C}$. Then: \textup{(i)} $L_\Omega(\eta; u) \geq 0$, with equality if and only if $u = \hat{u}_\Omega(\eta)$;
\textup{(ii)} $\nabla_{\eta} L_\Omega(\eta; u) = \hat{u}_\Omega(\eta) - u$;
\textup{(iii)} $\eta \mapsto L_\Omega(\eta; u)$ is convex.
\end{proposition}
\begin{proof}
\textup{(i)} The Fenchel--Young inequality $\Omega^*(\eta) + \Omega(u) \geq \eta^\top u$ is tight iff $\eta \in
\partial\Omega(u)$, which by the subgradient characterization of the conjugate is equivalent to $u = \nabla\Omega^*(\eta) = \hat{u}_\Omega(\eta)$
\citep{blondel2020learning}.
\textup{(ii)} $\nabla_{\eta} \Omega^*(\eta) = \hat{u}_\Omega(\eta)$ by Danskin's theorem and $-\eta^\top u$ is linear.
\textup{(iii)} $\Omega^*$ is convex as a supremum of affine functions; the rest is affine in $\eta$.
\end{proof}

\subsection{Perturbation and its equivalence with regularization}

For $\varepsilon>0$ and a noise distribution $\nu$ satisfying the conditions of Section~\ref{subsec:fy}, let $F_\varepsilon(\eta) := \E_{Z\sim\nu}[\max_{u \in \mathcal{Y}}(\eta+\varepsilon Z)^\top u]$ be the perturbed value function, so that $\hat{u}_\varepsilon(\eta)$ in~\eqref{eq:perturbed} is the expected perturbed maximizer. Based on Propositions~2.1 and~2.2 and Definition~4.1 of \citet{berthet2020learning} and Proposition~3.2 of \citet{dalle2022learning}, the following theorem states the equivalence between Fenchel--Young regularization and perturbation.

\begin{theorem}[\citealt{dalle2022learning}]
\label{thm:equiv} Let $\mathcal{Y}$ be compact, let $\nu$ satisfy the conditions of Section~\ref{subsec:fy}, and let $\varepsilon>0$. Set $F_1(\eta) := \E_{Z\sim\nu}[\max_{u\in\mathcal{Y}}(\eta+Z)^\top u]$ and let $\Omega := F_1^*$ be the regularizer induced by $\nu$, with domain $\mathcal{C'} := \dom(\Omega) \subseteq \mathcal{C}$. Then
\begin{equation}
  F_\varepsilon(\eta) = \sup_{u\in\mathcal{C'}}
    \bigl\{\eta^\top u - \varepsilon\,\Omega(u)\bigr\},
  \qquad
  \hat{u}_\varepsilon(\eta) = \arg\sup_{u\in\mathcal{C'}}
    \bigl\{\eta^\top u - \varepsilon\,\Omega(u)\bigr\}.
  \label{eq:equiv}
\end{equation}
\end{theorem}

In our application, the feasible set is $\mathcal{Y}=\mathcal{Y}_t$, which depends on the availability $a^t$ and, in the budgeted variant, on the budget $B^t$ (see~\eqref{eq:Yfeas_budget}); the induced regularizer depends on them as well, and the theorem applies at each stage with the corresponding $\mathcal{Y}_t$. By Theorem~\ref{thm:equiv}, the perturbed-optimizer gradient $\hat{u}_\varepsilon(\eta) - u$ is the gradient of the FY loss $L_{\varepsilon\Omega}(\eta; u) = F_\varepsilon(\eta) + \varepsilon\Omega(u) - \eta^\top u$, so training with perturbed optimizers minimizes this loss.

The proof below, reproduced in simplified form from the convex-analytic results of \citet{berthet2020learning} and \citet{dalle2022learning}, holds under the stated density and moment assumptions. It allows any compact feasible set; in the unbudgeted application the convex hull is the face $\{r \in [0,1]^{N_G+N_E} : r \le \mathbf{1}-a^t\}$ of the hypercube, and statements about interior points refer to the relative interior of that face.

\begin{proof}[Proof of Theorem~\ref{thm:equiv}]
$_{}$

\textit{Step~1 (convexity).} For each fixed $z \in \R^m$, the map $\eta \mapsto \max_{u\in\mathcal{Y}}(\eta+\varepsilon z)^\top u$ is the support function of $\mathcal{Y}$ evaluated at $\eta+\varepsilon z$, hence a supremum of affine functions of $\eta$ indexed by $u\in\mathcal{Y}$, and therefore convex; no finiteness of $\mathcal{Y}$ is assumed, and compactness is used only to ensure that the supremum is attained. Since $F_\varepsilon$ is an expectation of convex functions, $F_\varepsilon$ is convex.

\textit{Step~2 (gradient).} $\nabla F_\varepsilon(\eta) = \hat{u}_\varepsilon(\eta)$ by Danskin's theorem under the density and moment conditions of Section~\ref{subsec:fy}, which make the maximizer almost surely unique \citep{berthet2020learning}.

\textit{Step~3 (scaling).} $F_\varepsilon(\eta) = \varepsilon\,F_1(\eta/\varepsilon)$, obtained by factoring $\varepsilon$ inside the expectation.

\textit{Step~4 (conjugate).} $(F_\varepsilon)^*(u) = \sup_{\eta}\{\eta^\top u - \varepsilon F_1(\eta/\varepsilon)\} = \varepsilon\sup_v\{v^\top u - F_1(v)\} = \varepsilon\,\Omega(u)$, with the substitution $v = \eta/\varepsilon$.

\textit{Step~5 (Fenchel--Moreau).} $F_\varepsilon$ is convex (Step~1) and finite-valued on all of $\R^m$ (hence proper and lower semicontinuous). The Fenchel--Moreau theorem therefore gives $F_\varepsilon = (F_\varepsilon)^{**} = (\varepsilon\,\Omega)^*$, giving the first identity in~\eqref{eq:equiv}.

\textit{Step~6 (uniqueness).} $F_\varepsilon$ is differentiable (Step~2), so the supremum in Step~5 is uniquely attained. By standard convex analysis, the unique maximizer of $u \mapsto \eta^\top u - \varepsilon\,\Omega(u)$ equals the gradient of the conjugate $(\varepsilon\,\Omega)^*=F_\varepsilon$ at $\eta$, which by Step~2 is $\hat{u}_\varepsilon(\eta)$.
\end{proof}

\subsection{Proofs for the main paper}
\label{ec:proofs}

\begin{proof}[Proof of Proposition~\ref{prop:conditional}]
Fix a stage $t$ and a reachable pair $(x_t,\xi_t)$. Under the premise, the objective of the inner minimization in~\eqref{eq:bellman_decoupled} equals
\[
  Q_t^0(x_t,u_t,\xi_t) + \beta_t^\top u_t + h_t(x_t,\xi_t) - \bmu_t^*(x_t,\xi_t)^\top u_t ,
\]
which differs from the objective of the surrogate~\eqref{eq:surrogate_u} at $\bmu_t = \pi_{\bpsi}(x_t,\xi_t) = \bmu_t^*(x_t,\xi_t)$ only by $h_t(x_t,\xi_t)$, a term independent of $u_t$. Hence every minimizer of the surrogate attains the minimum in~\eqref{eq:bellman_decoupled}, and, together with an inner decision $\delta_t^0$ attaining $Q_t^0$ in~\eqref{eq:inner_cost}, the DFL decision attains the minimum in the Bellman equation~\eqref{eq:bellman}. The DFL policy therefore selects a Bellman-optimal decision at every stage and every reachable state, and by the standard verification argument for~\eqref{eq:bellman} \citep{shapiro2009lectures} its expected cost equals the optimal value of~\eqref{eq:msp}.
\end{proof}

\begin{proof}[Proof of the threshold rule~\eqref{eq:threshold}]
\textit{Feasibility.} If $a^t[k]=1$ then $r^t[k]\leq 1-a^t[k]=0$ forces $r^t[k]=0$. If $a^t[k]=0$ then $r^t[k]\in\{0,1\}$ is free.

\textit{Separability.} $(\beta-\bmu^t)^\top r^t = \sum_k (\beta[k]-\bmu^t[k])\,r^t[k]$, so the problem separates over components. For failed components ($a^t[k]=0$), minimize $(\beta[k]-\bmu^t[k])\,r^t[k]$ over $r^t[k]\in\{0,1\}$:
\begin{itemize}[nosep]
\item If $\bmu^t[k]>\beta[k]$: coefficient $<0$, so $r^t[k]=1$ is optimal.
\item If $\bmu^t[k]<\beta[k]$: coefficient $>0$, so $r^t[k]=0$ is optimal.
\item If $\bmu^t[k]=\beta[k]$: both are optimal; convention sets $r^t[k]=0$.
\end{itemize}
Combining over all $k$ yields~\eqref{eq:threshold}.
\end{proof}

\paragraph{Slack budgets.} Consider the budgeted surrogate~\eqref{eq:surr_knap}. If $B^t\ge\sum_{k:\,\bmu^t[k]>\beta[k],\ a^t[k]=0}\beta[k]$, the budget is inactive at the nominal values, and the deployment optimum coincides with the threshold rule~\eqref{eq:threshold}. Recovering a componentwise closed-form smoothing (with the Gaussian CDF $\Phi$ in place of $\sigma$ in~\eqref{eq:smooth_threshold}, since~\eqref{eq:perturbed_knap} uses Gaussian noise) requires the stronger condition $B^t\ge\sum_{k:\,a^t[k]=0}\beta[k]$, because a perturbed value can make a nominally unprofitable component profitable.

\section{Learning Implementation}
\label{ec:impl}

This section details the training settings and the relationship between the mathematical formulation of the main paper and the accompanying implementation, available at
\url{https://github.com/DavidKrame/DFL_for_MSP__Application_to_power_network_repair}.

\subsection{Training settings}

The network has three hidden layers of widths 128, 256, and 128 with ReLU activations and dropout 0.3, followed by a linear output layer. It is trained with Adam~\citep{kingma2014adam} (learning rate $10^{-3}$, weight decay $10^{-5}$, mini-batches of 64, gradient-norm clipping at 1), a ReduceLROnPlateau scheduler (patience 200, factor 0.5), and early stopping (patience 330 epochs). The maximum number of epochs scales with the training size, $\max(1,\lfloor n/10\rfloor)$ (200 epochs at $n=2{,}000$), so early stopping is active only for the larger sizes. The perturbation temperature is $\varepsilon=0.1$. The stage-level examples of the $n$ training scenarios are split at random into 80\%~training and 20\%~validation with a fixed seed; the validation value of the loss~\eqref{eq:fy_pseudo} drives the scheduler and early stopping, and the weights with the lowest validation value are restored at the end. Training scenarios are drawn with seed 42 and the 500 evaluation paths with seed 1000. The five runs per size use initialization seeds 1, 2, 3, 4, and 42; the data, the split, the dropout masks, and the batch order are the same in all five. Training labels come from single-scenario MIPs solved to a relative gap of $10^{-6}$, and in the budgeted variant the perturbed knapsack oracle~\eqref{eq:perturbed_knap} is estimated with $M=5$ Monte Carlo samples per training example.

\subsection{Perturbation distributions}
\label{appendix:impl:perturbation}

The two variants use different noise distributions. Without the budget, the surrogate~\eqref{eq:surrogate_repair} is separable, so the perturbed optimizer~\eqref{eq:perturbed} factorizes component-wise. For each failed component $k$ the inner problem is
\[
  \arg\max_{r_k \in \{0,1\}}
  \bigl(\eta[k] + \varepsilon\,z_k\bigr)\,r_k
  \;=\; \mathbf{1}\!\bigl[\eta[k] + \varepsilon\,z_k > 0\bigr],
\] and its expected value is $P(z_k > -\eta[k]/\varepsilon)$:
\begin{itemize}[nosep]
\item \textbf{Logistic noise} $z_k \sim \mathrm{Logistic}(0,1)$: $P(z_k > -\eta[k]/\varepsilon) = \sigma(\eta[k]/\varepsilon)$, the smoothing~\eqref{eq:smooth_threshold} used in our unbudgeted runs. A standard logistic variable is the difference of two independent standard Gumbel variables, so the binary repair choice takes the Gumbel-max form, with the sigmoid as its choice probability.
\item \textbf{Gaussian noise} $z_k \sim \mathcal{N}(0,1)$: $P(z_k > -\eta[k]/\varepsilon) = \Phi(\eta[k]/\varepsilon)$, with $\Phi$ the standard normal CDF.
\end{itemize}
With the budget, the knapsack oracle is perturbed with standard Gaussian noise, as in~\eqref{eq:perturbed_knap}, and the expectation is estimated by Monte Carlo. Both choices are valid instantiations of the perturbed-optimizer framework, and Theorem~\ref{thm:equiv} holds for either, with the corresponding regularizer $\Omega = F_1^*$ determined by the noise distribution. The gradient signal $\hat{r}_\varepsilon(\eta) - r^*$ has the same sign structure in both cases.

\subsection{Stop-gradient surrogate for the FY loss}
\label{appendix:impl:pseudoloss}

Stochastic gradient descent uses only the gradient of the loss, and that gradient, $\hat{r}_\varepsilon(\eta) - r^*$, needs neither $\Omega$ nor the value of $L_{\varepsilon\Omega}$. With the budget, $\Omega$ is the regularizer induced by the Gaussian noise, which has no closed form, and the knapsack solutions behind $\hat{r}_\varepsilon$ are piecewise constant in $\eta$, so automatic differentiation has nothing to differentiate through. In place of the FY loss $L_{\varepsilon\Omega}$ in~\eqref{eq:train_loss}, we therefore train on the simplest scalar whose automatic gradient is the FY gradient, with one implementation for both variants:
\begin{equation}
  L_{\mathrm{pseudo}}(\eta;\, r^*)
  \;:=\;
  \bigl(\hat{r}_\varepsilon(\eta).\mathrm{detach}()
  - r^*\bigr)^\top \eta,
  \label{eq:fy_pseudo}
\end{equation}
where \texttt{.detach()} treats the smoothed decision as a constant (stop-gradient). Its gradient,
\begin{equation}
  \nabla_\eta\, L_{\mathrm{pseudo}}(\eta;\, r^*)
  = \hat{r}_\varepsilon(\eta) - r^*,
  \label{eq:pseudo_grad}
\end{equation}
coincides with $\nabla_\eta\, L_{\varepsilon\Omega}$ from Proposition~\ref{prop:fy_props}(ii) and Theorem~\ref{thm:equiv}, so a gradient step on $L_{\mathrm{pseudo}}$ coincides with a gradient step on $L_{\varepsilon\Omega}$ at the same step size: the estimator is the perturbed-optimizer gradient of the FY loss, neither a straight-through construction nor the SPO\texttt{+} surrogate. The loss values differ, however. $L_{\varepsilon\Omega}$ is nonnegative and vanishes only when the smoothed decision equals the target (Proposition~\ref{prop:fy_props}(i)), whereas $L_{\mathrm{pseudo}}$ can be negative and can vanish when $\hat{r}_\varepsilon \neq r^*$ (if $\hat{r}_\varepsilon - r^*$ is orthogonal to $\eta$). The training curves are therefore pseudo-loss values, and because the learning-rate schedule, early stopping, and checkpoint choice use the validation pseudo-loss, the training trajectory can differ from one driven by the theoretical FY loss itself.

In the separable, unbudgeted case with the logistic smoothing~\eqref{eq:smooth_threshold}, the update $\sigma(\eta[k]/\varepsilon) - r^*[k]$ of a failed component $k$ is $\varepsilon$ times the gradient in $\eta[k]$ of the binary cross-entropy between $r^*[k]$ and $\sigma(\eta[k]/\varepsilon)$, so it coincides, up to a constant absorbed into the learning rate, with cross-entropy training on the same labels. The coincidence concerns the update rule only, since $\sigma(\eta[k]/\varepsilon)$ is a smoothed oracle decision rather than a classifier's probability, and it disappears under the budget, where $\hat{r}_\varepsilon$ no longer factorizes across components.

\section{Full Model and Failure Process}
\label{ec:model}

\subsection{The complete multistage program}

Combining the stage DC-OPF~\eqref{eq:dc_opf} and the dynamics~\eqref{eq:transition}, the repair-scheduling problem~\eqref{eq:msmip} written out in full is
\begin{equation}
\begin{aligned}
  \min_\pi \quad
  & \E\!\left[\sum_{t=1}^T \bigl(c^\top p_g^t + \alpha\,\mathbf{1}^\top(p_l^t - d^t)
    + \beta^\top r^t\bigr)\right] \\
  \text{s.t.}\quad
  & \text{\eqref{eq:transition} holds }\forall\,t,
    \text{ with }a^0=\mathbf{1},\;r^0=\mathbf{0}, \\
  & \theta^t[\mathrm{ref}] = 0,\quad
    -\theta^{\max}\mathbf{1} \leq \theta^t \leq \theta^{\max}\mathbf{1}, \;\forall\,t, \\
  & C f^t = S_g p_g^t - S_l d^t, \;\forall\,t, \\
  & -(\mathbf{1}-a_e^t) \odot M \leq f^t - BC^\top\theta^t
    \leq (\mathbf{1}-a_e^t) \odot M, \;\forall\,t, \\
  & -a_e^t \odot f^{\max} \leq f^t
    \leq a_e^t \odot f^{\max}, \;\forall\,t, \\
  & a_g^t \odot \underline{p}_g \leq p_g^t \leq a_g^t \odot \overline{p}_g, \;\forall\,t, \\
  & \mathbf{0} \leq d^t \leq p_l^t,\;\forall\,t, \\
  & \mathbf{0} \leq r^t \leq \mathbf{1}-a^t,\;\;
    r^t \in \{0,1\}^{N_G+N_E},\;\forall\,t,
\end{aligned}
  \label{eq:msmip_full}
\end{equation}
with each $d^t$ and $r^t$ $\sigma(\xi^1,\ldots,\xi^t)$-measurable. The constraints are those of~\eqref{eq:dc_opf} with symmetric line-flow capacities $f^{\max}\in\mathbb{R}^{N_E}_{>0}$ ($\overline{f} = f^{\max}$, $\underline{f} = -f^{\max}$); the budgeted variant adds~\eqref{eq:budget} at every stage.

\subsection{Failure probability model}

In~\eqref{eq:logit}, the parameters $A[k]$ and $B[k]$ are determined by two component inputs, the peak failure stage $t^*[k]$ and the peak probability $p_\text{max}[k]$, via
\begin{equation}
  B[k] = \frac{\logit(p_\text{max}[k])}{\delta[k]},
  \qquad
  A[k] = \frac{t^*[k]+K}{\delta[k]}\,\logit(p_\text{max}[k]),
  \label{eq:logit_params}
\end{equation}
where $\delta[k] := (t^*[k]+K)\ln(t^*[k]+K) - t^*[k]$. These expressions are obtained by requiring $\logit(\lambda^{t^*[k]}[k]) = \logit(p_\text{max}[k])$ and $\partial\logit(\lambda^t[k])/\partial t\big|_{t=t^*[k]} = 0$. They are well defined when $\delta[k] > 0$; since $\delta[k]$ increases with $K$, for $t^*[k]=1$, as in our experiments, this holds if and only if $K > K^*\approx 0.763$, and hence for the experimental value $K=1.5$. In addition, the condition $p_\text{max}[k] > \tfrac{1}{2}$ ensures $B[k] > 0$, so that $\lambda^t[k] \to 0$ as $t \to \infty$.

\begin{figure}[tb]
{\centering
\includegraphics[width=0.82\linewidth]{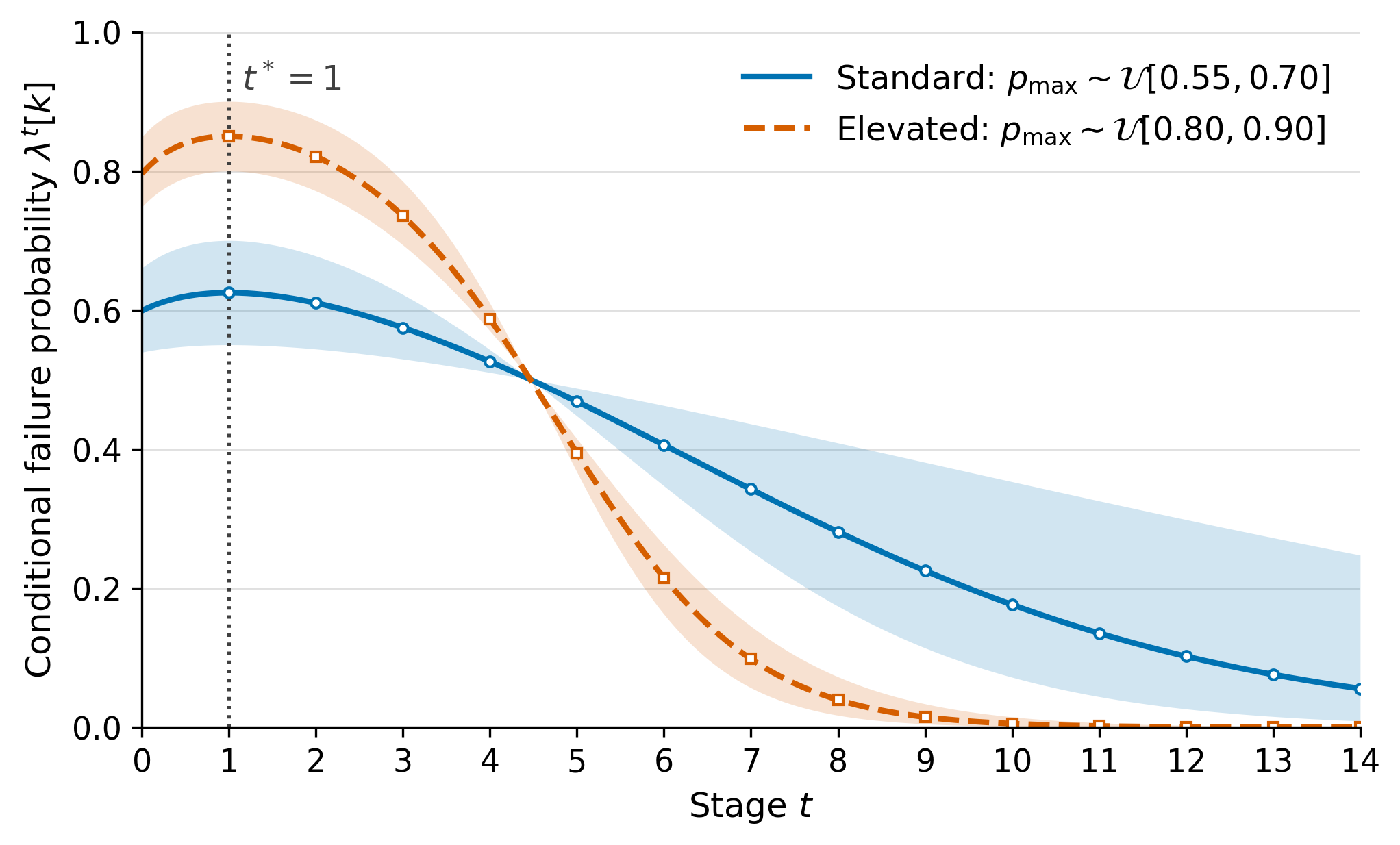}\par}
\caption{Experimental failure-probability profiles.\label{fig:failure_curve}}
{\small\textit{Note.} Profiles use $T=14$, $t^*=1$, and $K=1.5$. Each component is assigned with probability $1/2$ to the standard class ($p_{\max}\sim\mathcal U[0.55,0.70]$) or elevated class ($p_{\max}\sim\mathcal U[0.80,0.90]$). Shading spans each class's $p_{\max}$ range; marked curves use its midpoint. Probabilities are conditional on the component being operational. All profiles cross $1/2$ between stages 4 and 5 (at $t\approx4.46$), whatever their peak value; a larger $p_{\max}$ gives a higher peak and a faster decay afterward.\par}

\end{figure}

Substituting~\eqref{eq:logit_params} into~\eqref{eq:logit} gives
\[
  \logit(\lambda^t[k]) = \frac{\logit(p_\text{max}[k])}{\delta[k]}\,\bigl[(t^*[k]+K)\ln(t+K) - t\bigr],
\]
whose sign does not depend on $p_\text{max}[k]$ when $p_\text{max}[k] > \tfrac12$. With $t^*=1$ and $K=1.5$, the bracket is positive for $t\le4$ and negative for $t\ge5$ (it vanishes at $t\approx4.46$), so every failure probability is above $1/2$ up to stage 4 and below it from stage 5 on (Figure~\ref{fig:failure_curve}).

\section{Comparators and Classical Methods}
\label{ec:comparators}

\subsection{SAA and SDDiP}

\emph{SAA} solves a single extensive-form MILP on a sampled binary scenario tree with branching factor $q=2$ over $T=14$ stages, giving $2^{14}=16{,}384$ leaf paths in every configuration. The MILP is solved to a target relative gap of $10\%$ or until the solver time limit just below the shared wall-clock limit. For evaluation on the common paths, the resulting node policy is rolled out stage by stage: at each stage, the realized path is mapped to the current node's closest child in Hamming distance, and that child's repair decision is applied subject to the realized component availability.

\emph{SDDiP} runs Lagrangian cut generation through \texttt{SDDP.jl}~\citep{dowson2021sddpjl} using HiGHS and terminates at the shared wall-clock limit, without gap-based stopping. We report the final lower bound and the final policy's mean cost over 500 post-training simulations drawn from SDDiP's stagewise distribution, which enumerates the full failure law on 5-bus and uses a sampled finite support on case14. Their difference is an internal estimated policy gap, which may favor SDDiP; because the runs did not converge within the shared budget, it is the most informative assessment of the final policies that these runs provide.

\begin{table}[tb]
\caption{SDDiP results at termination of the wall-clock budget.\label{tab:sddip}}
{\centering
\begin{tabular}{@{}llrrrr@{}}
\toprule
Network & Variant & Iterations & LB (\$) & Policy mean (\$) & Est.\ gap \\
\midrule
5-bus ($s{=}100$)   & unbudgeted & 782 & 1{,}024{,}482 & 1{,}584{,}476 & 35.3\% \\
5-bus ($s{=}100$)   & budgeted   & 755 & 1{,}296{,}240 & 1{,}651{,}543 & 21.5\% \\
case14 ($s{=}45$) & unbudgeted & 158 &   671{,}289 & 2{,}110{,}204 & 68.2\% \\
case14 ($s{=}45$) & budgeted   & 150 &   831{,}741 & 2{,}208{,}569 & 62.3\% \\
case14 ($s{=}55$) & unbudgeted & 150 &   740{,}410 & 2{,}283{,}357 & 67.6\% \\
case14 ($s{=}55$) & budgeted   & 147 &   927{,}540 & 2{,}304{,}480 & 59.8\% \\
\bottomrule
\end{tabular}\par}
\smallskip
{\small\textit{Note.} LB is the final SDDiP lower bound. The policy mean is reconstructed from the saved stagewise averages over 500 post-training simulations, and the estimated gap is $(\text{policy mean}-\text{LB})/\text{policy mean}$. The simulations use SDDiP's own stagewise distribution, which is exact on 5-bus and sampled on case14, rather than the common evaluation paths used in Table~\ref{tab:dfl}. Individual simulation costs were not retained, so no confidence intervals are reported.\par}

\end{table}

\subsection{The near-sighted MPC baseline}
\label{appendix:mpc}

This section specifies the near-sighted model-predictive control (MPC) baseline of Section~\ref{subsec:setup}, a non-anticipative, learning-free controller. Standard stochastic MPC repeatedly solves a finite-horizon expected-cost problem with a terminal value, implements the first action, and replans \citep{bertsekas2024MPC}; our baseline keeps this receding-horizon feedback but replaces the expectation by a mean-failure rollout and sets the terminal value to zero.

\paragraph{Setup.} We use the notation of Section~\ref{subsec:msmip}, write $Q^0(a)$ for the optimal value of the DC-OPF~\eqref{eq:dc_opf} at availability $a$ (Proposition~\ref{prop:decoupling}), and write $\beta(S)=\sum_{k\in S}\beta[k]$ for a repair set $S$. The rollout below evaluates $Q^0$ at fractional $a\in[0,1]^{N_G+N_E}$ on a derated network, in which $a_e[j]$ scales both the susceptance and the thermal capacity of line $j$, $f[j]=a_e[j]\,b_j\bigl(\theta[\mathrm{from}(j)]-\theta[\mathrm{to}(j)]\bigr)$ and $|f[j]|\le a_e[j]\,f^{\max}[j]$, and $a_g[i]$ scales the capacity of generator $i$. This heuristic, used only to score schedules, is a linear program for fixed $a$ and coincides with~\eqref{eq:dc_opf} at binary $a$.

\paragraph{Our MPC baseline.} The baseline is parameterized by the window $W$ with $\tau+1\le W\le T$ and by a cap $N_{\max}$ on the number of candidate schedules scored per stage; the experiments use $W=3$ and $N_{\max}=2187=3^7$. At each true stage $t$ it performs the following.
\begin{enumerate}[nosep]
\item \emph{Observe.} Read $a^t$ and the currently-failed set $F^t=\{k:a^t[k]=0\}$; only these are repair candidates.
\item \emph{Effective window.} Let $L_t=\min(W,\,T-t+1)$. A repair ordered at window-step $j$ completes at $j+\tau$, so it can influence the window only if $j+\tau\le L_t-1$; the \emph{actionable steps} are $j\in\{0,\dots,L_t-1-\tau\}$. If none exist (final stages) or $F^t=\varnothing$, order no repair.
\item \emph{Build the candidate set.} A schedule assigns each failed component to an actionable step or to ``never'': $\varphi:F^t\to\{0,\dots,L_t-1-\tau\}\cup\{\infty\}$, with $S_j=\varphi^{-1}(j)$. There are $(L_t-\tau+1)^{|F^t|}$ such schedules. If this number does not exceed $N_{\max}$, the candidate set $\mathcal A_t$ is the full set of schedules, enumerated in lexicographic order from ``repair every failed component now'' to the all-never schedule $\varphi\equiv\infty$. Otherwise $\mathcal A_t$ holds the all-never schedule, listed first, and $N_{\max}-1$ further schedules, each drawn uniformly from the full set (with replacement) by a random-number generator seeded once with seed $0$.
\item \emph{Score under expected failures.} Collapse the future to its mean: set $\bar a^{(0)}=a^t$ and, for $j=0,\dots,L_t-2$,
\begin{equation}
  \bar a^{(j+1)}[k]=\bar a^{(j)}[k]\bigl(1-\lambda^{t+j+1}[k]\bigr)
  +\bigl(1-\bar a^{(j)}[k]\bigr)\,\mathbf 1\!\left[k\in S_{j+1-\tau}\right],
  \label{eq:ce_rollout}
\end{equation}
i.e.\ the dynamics~\eqref{eq:transition} with $\xi$ replaced by its mean $\lambda$. Availabilities become fractional, which $Q^0(\cdot)$ accepts through the derated network defined above. The schedule's predicted window cost is
\begin{equation}
  \widehat J_t(\varphi)=\sum_{j=0}^{L_t-1}\Bigl[Q^0\bigl(\bar a^{(j)}\bigr)+\beta(S_j)\Bigr].
  \label{eq:mpc_windowcost}
\end{equation}
\item \emph{Optimize over the candidate set, execute only the current stage.} Choose $\varphi^\star$ as the first minimizer of $\widehat J_t$ over $\mathcal A_t$ in the order above: an exact minimization over all schedules when the cap does not bind, and a randomized search over the sampled subset when it does. The all-never schedule belongs to $\mathcal A_t$ by construction, so the chosen schedule never has a higher predicted window cost~\eqref{eq:mpc_windowcost} than ordering no repair. Ties are resolved by this order: in favor of earlier repairs under exact enumeration, and in favor of the all-never schedule in the sampled search. Execute only the current-stage repairs $r^t=\mathbf 1[\,\cdot\in S_0^\star\,]$ and discard $S_1^\star,\dots,S_{L_t-1}^\star$.
\item \emph{Advance and replan.} Apply $r^t$, observe the realized $\xi^{t+1}$, update $a^{t+1}$ via~\eqref{eq:transition}, and return to step~1 at $t+1$. Replanning from the observed state corrects for the difference between the mean rollout and the realized failures.
\end{enumerate}

\paragraph{Scope and variants.} Relative to standard stochastic MPC, the baseline has zero terminal value, uses the certainty-equivalent rollout~\eqref{eq:ce_rollout}, and schedules only components failed at the current stage. It therefore has no preemptive repairs and credits no benefits beyond the window. The search is exhaustive only while $(L_t-\tau+1)^{|F^t|}\le N_{\max}$; with $W=3$, $\tau=1$, and $N_{\max}=3^7$, this covers at most seven simultaneously failed components when the full window is available (more in the final stages, where the window is shorter), after which the fixed-seed sample is used. Under the per-stage budget~\eqref{eq:budget}, step~3 retains only schedules satisfying $\beta(S_j)\le B^{t+j}$ at every window step; the empty schedule preserves feasibility, and scoring and replanning are unchanged.

\section{Computing Resources}
\label{ec:compute}

All runs used the nodes and the per-job wall-clock limit described in Section~\ref{subsec:setup}. The offline methods differ in how the 64 cores are used: DFL label generation runs up to 64 independent single-threaded MILP solves in parallel (HiGHS), DFL training is a single CPU process, SAA runs one multi-threaded extensive-form solve (Gurobi, which used up to 32 of the cores), and SDDiP gives the cores to its stage-subproblem solver (HiGHS). The stage DC-OPF problems are solved with HiGHS, except in the SAA rollouts, which use Gurobi. Near-sighted MPC has no offline stage and is timed online only, at deployment. The limit applies to each job, so it does not equalize total effort across methods: DFL runs separate jobs for label generation and for each training size and seed.

\end{document}